\documentclass[a4paper,10pt]{article}

 \usepackage{amsmath, amssymb, latexsym, amscd, amsthm,amsfonts,amstext}
 \usepackage[mathscr]{eucal}
 \usepackage{graphicx}
 \usepackage{subfig}

\title{Globally convergent and adaptive finite element methods  in imaging of buried objects from experimental backscattering radar measurements}

\author{Larisa Beilina$^{1}$\footnote{Corresponding author}, Nguyen Trung Th\`anh$^{2}$, Michael V. Klibanov$^{3}$ \\
and John Bondestam Malmberg$^{1}$\\
$^{1}$Department of Mathematical Sciences, \\
Chalmers University of Technology and
Gothenburg University, \\
SE-42196 Gothenburg, Sweden \\
Emails: \texttt{larisa@chalmers.se, john.bondestam.malmberg@chalmers.se}.\\
$^2$Department of Mathematics, Iowa State University, \\
 Ames, IA 50011, USA. Email: \texttt{thanh@iastate.edu}.\\
$^3$Department of Mathematics \& Statistics, University of North Carolina at Charlotte,\\
Charlotte, NC 28223, USA. Email: \texttt{mklibanv@uncc.edu}.}
\date{}

\begin{document}
\maketitle

\begin{abstract}

We consider a two-stage numerical procedure for imaging of objects buried
in dry sand using time-dependent backscattering experimental
radar measurements. These measurements are generated by a single point
source of electric pulses and are collected using a microwave
scattering facility which was built at the University of North
Carolina at Charlotte.
Our imaging problem is formulated as the inverse problem of the
reconstruction of the spatially distributed dielectric permittivity
$\varepsilon _\mathrm{r}\left( \mathbf{x}\right), \ \mathbf{x}\in
\mathbb{R}^{3}$, which is an unknown coefficient in Maxwell's
equations.

 On the first stage an approximately globally convergent method
  is applied to get a good first approximation for the exact
 solution.  On
 the second stage a local adaptive finite element method is applied to refine the solution obtained on the
 first stage. The two-stage numerical procedure results in accurate
 imaging of all three components of interest of targets: shapes,
 locations and refractive indices. In this paper we briefly describe
 methods and present new reconstruction results for both stages.

\end{abstract}

\textbf{Keywords}: Inverse scattering, refractive indices, approximately globally convergent algorithm, adaptive finite element method.

\textbf{AMS classification codes:} 65N15, 65N30, 35J25.

\graphicspath{{FIGURES/}
 {pics/}}

\section{Introduction}
\label{sec:1}

In this paper we consider the problem of reconstruction of refractive
indices, shapes and locations of buried objects in the dry sand from
backscattering time-dependent experimental data using the two-stage
numerical procedure presented in \cite{BKK, BK, BTKJ, BK_AA}. Our problem is a
coefficient inverse problem (CIP) for  Maxwell's
equations in three dimensions. Experimental data were collected using a microwave
scattering facility which was built at the University of North
Carolina at Charlotte, USA.  Our experimental data are collected using
a single location of the source. The backscattered signal is measured
on a part of a plane. Our potential applications are in
the imaging of explosives, such as land mines and improvised
explosive devices.
This work is a continuation of our recent works on this topic, where
we have treated a much simpler case of experimental data for targets
placed in air \cite{BTKF, BTKJ, NBKF}.

The two-stage numerical procedure means that we combine two different
methods to solve our CIP.  On the first stage the approximately
globally convergent method of \cite{BK} is applied in order to obtain
a good first approximation for the exact solution.  We have presented
results of reconstruction on this stage in our publications
\cite{BTKF, NBKF} for objects placed in air. 
In our recent study \cite{NBKF2} we presented reconstructions of twenty-five
objects which show that the method of \cite{BK} works well in
estimating the dielectric constants (equivalently, refractive indices)
and locations of buried objects.

 In \cite{KBB} it was investigated why a minimizer of the Tikhonov
 functional is indeed closer to the exact solution than the first
 guess of this minimizer. Because of that it makes sense improve the
 solution which we have obtained on the first stage of our two-stage
 numerical procedure. To do that  the local adaptive finite element method of
 \cite{BMaxwell2} is applied by taking the solution of the first stage
 as the starting point in the minimization of a Tikhonov functional in
 order to obtain better approximations and shapes of objects on the
 adaptively refined meshes. In \cite{BTKJ} it was shown that using the
 adaptive finite element method all three components of interest for
 targets placed in the air can be simultaneously imaged: refractive
 indices, shapes and locations.

 Compared to imaging of targets placed in the air (see \cite{BTKF, BTKJ,
   NBKF}), there are three main difficulties in imaging of buried
 targets: (i) the  signals of targets are much weaker than those when the
 targets are in air, (ii) these signals may overlap with the
 reflection from the ground's surface, which makes them difficult to
 distinguish, and (iii) the reflection from the grounds surface may
 dominate the target's signals after the Laplace transform since the
 kernel of the Laplace transform decays exponentially with respect to
 time. We have handled this difficulty in \cite{NBKF2} via a new data
 preprocessing procedure. This procedure results in preprocessed
 data, which are used as the input for our globally convergent
 algorithm, that is, the input for the first stage of our method.

It is notable that we have experimentally observed a rare superresolution phenomenon and have numerically reconstructed the
corresponding image (see section \ref{sec:numex}). The resolution
limit which follows from the Born approximation, that is, in the
diffraction limit, is $\lambda/2$, where $\lambda$ is the wavelength
of the signal. However, we have resolved two targets with the distance
$\lambda/4.5$ between their surfaces. It was shown in, for instance
\cite{Simonetti:APL2006}, that the superresolution can occur because of
nonlinear scattering, and our algorithm is nonlinear, including the
step of extraction of the target's signal in our data preprocessing
procedure  \cite{NBKF2}. Experimentally the superresolution
phenomenon was demonstrated in \cite{CC}. We also refer to the recent
work \cite{Ammari} where the superresolution is discussed.

An outline of this paper follows. In section \ref{sec:stage1} we
briefly describe the approximate globally convergent method.  In section
\ref{sec:stage2} we present the forward, inverse, and adjoint problems as well as the Tikhonov functional for the second stage.  In section
\ref{sec:fem} we describe the finite element method used in
computations and in section \ref{sec:general} we investigate general
framework for a posteriori error estimation for CIPs.
 In section \ref{subsec:ad_alg} we describe the mesh refinement
 recommendation and the adaptive algorithm. In section \ref{sec:numex}
 we present results of our computations.

\section{The first stage}

\label{sec:stage1}

In this section we state the forward and inverse problems which we
consider on the first stage. We also briefly outline the globally
convergent method of \cite{BK} and present the algorithm used in
computations of the first stage.

\subsection{Forward and inverse problems}

Let $\Omega \subset \mathbb{R}^{3}$ be a convex bounded domain with
the boundary $\partial \Omega \in C^{3}.$ Denote the spatial coordinates by $\mathbf{x}=\left(
x,\,y,\,z\right) \in \mathbb{R}^{3}.$ Let $C^{k+\alpha }$ be
H\"older spaces, where $k\geq 0$ is an integer and $\alpha \in \left(
0,\,1\right).$ We consider the propagation of the electromagnetic wave
in $\mathbb{R}^{3}$ generated by an incident plane wave. On the first
stage we model the wave propagation by the following Cauchy problem
for the scalar wave equation
\begin{align}
&{\varepsilon_\mathrm{r}} (\mathbf{x})\frac{\partial^2u}{\partial t^2}(\mathbf{x},\,t) - \Delta u(\mathbf{x},\,t) = \delta(z- z_0) f(t),&&(\mathbf{x},\,t)\in \mathbb{R}^{3}\times (0,\,\infty ),  \label{eq:fp1} \\
&u(\mathbf{x},\,0)=0,\quad \frac{\partial u}{\partial t}(\mathbf{x},\,0)= 0, &&\mathbf{x}\in \mathbb{R}^3.  \label{eq:fp2}
\end{align}
 Here $f\left( t\right)\not\equiv 0$ is the time-dependent incident plane wave
 at the plane $\left\{
z=z_{0}\right\}$,  $u$ is the total wave generated by   $f(t)$  and
propagating along the $z$-axis.

Let the function $E\left( \mathbf{x},\,t\right) $ represent the voltage
of one component $E_2$ of the electric field $E\left(
\mathbf{x},\,t\right) =\left( E_1,\,E_2,\,E_3\right) \left(\mathbf{x},\,t\right).$
 In our experiments the component $E_2$
corresponds to the electromagnetic wave which is sent into the medium.
Our mathematical model of the first stage uses only the single
equation (\ref{eq:fp1}) with $u=E_2$ instead of the full Maxwell's system.  We can
do such approximation since  it was shown numerically in \cite{BM}  that
the component  $E_2$ of the electric field $E $ dominates the other two components in the case we consider. See also
\cite{BK} where a similar scalar wave equation was used to work with
transmitted experimental data.

The function ${\varepsilon_\mathrm{r}}$ in (\ref{eq:fp1}) represents the
spatially distributed dielectric permittivity. We assume that ${\varepsilon_\mathrm{r}}$ is unknown inside the domain $\Omega \subset\mathbb{R}^{3}$ and is such that
\begin{equation}
\varepsilon _\mathrm{r}\in C^{\alpha }\left( \mathbb{R}^{3}\right),\quad
\varepsilon _\mathrm{r}(\mathbf{x})\in \lbrack 1,\,b] \text{ for } \mathbf{x}\in \mathbb{R}^3,\quad
\varepsilon _\mathrm{r}(\mathbf{x})=1\text{ for }\mathbf{x}\in \mathbb{R}^{3}\setminus \Omega ,  \label{2.20}
\end{equation}
where $b>1$ is a constant. We assume that the set of
admissible coefficients in (\ref{2.20}) is known.  Let $\Gamma \subset \partial
\Omega $ be a part of the boundary $\partial \Omega.$ In our
experiments the plane wave is initialized outside of the domain
$\overline\Omega$, that is $\overline\Omega\cap \{z=z_{0}\} = \varnothing$.

\textbf{Coefficient Inverse Problem (CIP).}
\emph{ Determine the function }$\varepsilon _\mathrm{r}\left( \mathbf{x}\right) $\emph{ for }$\mathbf{x}\in \Omega $\emph{, assuming that the following function
}$g$\emph{ is known for a single incident plane wave
  generated at the plane} $\{z = z_{0}\}$ \emph{outside of} $\overline\Omega$:
\begin{equation}
u\left( \mathbf{x},\,t\right) = g\left( \mathbf{x},\,t\right)~ \forall \left(
\mathbf{x},\,t\right) \in \Gamma \times \left( 0,\,\infty \right) .\notag
\end{equation}

Global uniqueness theorems for multidimensional CIPs with a single
measurement are currently known only under the assumption that at
least one of initial conditions does not equal zero in the entire
domain $\overline{\Omega }$ \cite{BK, BukhK}.  However, this is not
our case and the method of Carleman estimates is inapplicable to our
CIP. Thus, we simply assume that uniqueness of our CIP holds.

\subsection{The globally convergent method}

\label{sec:gca}

Here we briefly present approximately globally convergent method of \cite{BK}.

We perform a Laplace
 transformation
\begin{equation}
\tilde u(\mathbf{x},\,s)=\int\limits_{0}^{\infty }u(\mathbf{x},\,t)e^{-st}\,\mathrm{d}t,
\notag
\end{equation}%
where $s$ is a positive parameter which we call \textit{pseudo
  frequency}. We assume that $s\geq \underline{s}>0$ and denote by
$\tilde f(s)$ the Laplace transform of $f(t)$.  We assume that $\tilde
f(s) \neq 0$ for all $s \ge \underline{s}$. Define $w(\mathbf{x},\,s):=
\tilde u(\mathbf{x},\,s)/\tilde f(s)$. The function $w$ satisfies the
equation
\begin{equation}
\Delta w(\mathbf{x},\,s)-s^{2}{\varepsilon_\mathrm{r}} (\mathbf{x})w(\mathbf{x},\,s)=-\delta (z - z_{0}),\quad\mathbf{x}\in \mathbb{R}^{3},\,\ s\geq
\underline{s}.  \label{eq:w}
\end{equation}
It was shown in \cite{NBKF2}
 that
$w(\mathbf{x},\,s)>0$ and $ \lim_{\left\vert \mathbf{x}\right\vert
  \rightarrow \infty }\left[ w\left( \mathbf{x},\,s\right)
  -w_0(\mathbf{x},\,s)\right] =0, $ where $w_{0}\left(
\mathbf{x},\,s\right) :=e^{-s\left\vert z-z_{0}\right\vert} /\left(
2s\right) $ is a solution of equation (\ref{eq:w}) for the case
${\varepsilon_\mathrm{r}} \left( \mathbf{x}\right) \equiv 1,$ which decays to zero as
$\left\vert z\right\vert \rightarrow \infty.$ Next, introduce the
function $v$ by $v(\mathbf{x},\,s):=\ln\big(w(\mathbf{x},\,s)\big)/s^{2}$  and substitute $w=\exp(vs^{2})$ into
(\ref{eq:w}). By noting  that $\overline\Omega \cap \{z =
z_0\} = \varnothing$, we obtain the following equation for the
explicit computation of the coefficient ${\varepsilon_\mathrm{r}}$:
\begin{equation}
\Delta v(\mathbf{x},\,s)+s^{2}|\nabla v(\mathbf{x},\,s)|^{2}={\varepsilon_\mathrm{r}} (\mathbf{x}),\quad\mathbf{x}\in \Omega,~ s\geq\underline{s} .
\label{eq:c}
\end{equation}
 Next, we eliminate the
unknown coefficient ${\varepsilon_\mathrm{r}} (\mathbf{x})$ from (\ref{eq:c}) by taking the
derivative with respect to $s$ both sides of (\ref{eq:c}). Denote by $
q:=\frac{\partial v}{\partial s}$, then
\begin{equation}
v(\mathbf{x},\,s)=-\int\limits_{s}^{\infty }q(\mathbf{x},\,\tau)\,\mathrm{d}\tau  = -\int\limits_{s}^{\bar{s}}q(\mathbf{x},\,\tau)\,\mathrm{d}\tau + V(\mathbf{x}),
\notag
\end{equation}
where $\bar{s}>\underline{s}$. We call the function
$V(\mathbf{x})=v(\mathbf{x},\,\bar{s})$ the \textquotedblleft
tail function\textquotedblright and define it  by
\begin{equation}
V(\mathbf{x})=\frac{\ln w(\mathbf{x},\,\bar{s})}{\bar{s}^{2}}.  \label{eq:tail}
\end{equation}
From (\ref{eq:c}) we obtain the following  equation  for two unknown
  functions $q$ and $V$
\begin{equation}
\begin{aligned}
\Delta q(\mathbf{x},\,s)& - 2s^{2}\nabla q(\mathbf{x},\,s)\cdot \int\limits_{s}^{\bar{s}}\nabla q(\mathbf{x},\,\tau )\,\mathrm{d}\tau + 2s^{2}\nabla V(\mathbf{x})\cdot \nabla q(\mathbf{x},\,s)\\
&+2s\left\vert \int\limits_{s}^{\bar{s}}\nabla q(\mathbf{x},\,\tau )\,\mathrm{d}\tau \right\vert ^{2} - 4s\nabla V(\mathbf{x})\cdot \int\limits_{s}^{\bar{s}}\nabla q(\mathbf{x},\,\tau)\,\mathrm{d}\tau + 2s\left\vert \nabla V(\mathbf{x})\right\vert ^{2}=0, \label{eq:q}
\end{aligned}
\end{equation}
for $\mathbf{x}\in\Omega$ and $s\in(\underline s,\,\bar s)$.

To find the tail function $V$ we use an iterative procedure presented
in the next section, see \cite{BTKF, NBKF} for details of this procedure.
The function $q$ satisfies the following boundary condition
\begin{equation}
q(\mathbf{x},\,s)=\psi (\mathbf{x},\,s),\quad\mathbf{x}\in \partial \Omega ,
\label{eq:q2}
\end{equation}
where
$
\psi (\mathbf{x},\,s)=\frac{\partial }{\partial s}\left[ \frac{\ln \varphi(\mathbf{x},\, s)}{
s^{2}}\right]$   with  $\varphi
(\mathbf{x},\,s)=\int\limits_{0}^{\infty }g(\mathbf{x},\,t)e^{-st}\,\mathrm{d}t /\tilde f(s)$.

\subsection{Iterative procedure and description of the approximate globally convergent algorithm}
\label{subsec:alg}

In our iterative procedure we
 divide the pseudo frequency interval $[\underline{s},\,\bar{s}]$ into
$N$ sub-intervals $\bar{s}=s_{0}>s_{1}>\cdots
>s_{N}=\underline{s}$  of the step size $h$  such that $s_{n}-s_{n+1}=h.$ We approximate the function
$q$ by a piecewise constant function with respect $s$,
$q(\mathbf{x},\,s)\approx q_{n}(\mathbf{x})$, $s\in
(s_{n},\,s_{n-1}],$ $n=1,\,\dots,\,N$,  and set $q_{0}\equiv 0$.
 Next,  we multiply equation \eqref{eq:q} by the  Carleman Weight
  Function $\exp \left[ \Lambda \left( s-s_{n-1}\right) \right]$, $s\in
  \left( s_{n},\,s_{n-1}\right)$, where $\Lambda \gg 1$ is a large
  parameter chosen in the computations,  and integrate with respect to $s$  over every  pseudo frequency
interval $[s_n,\, s_{n-1}]$.  Finally, we get
   a system of elliptic equations for the functions $q_{n}$ for $\mathbf{x}\in\Omega$:
\begin{equation}
\begin{aligned}
\Delta q_{n}(\mathbf{x}) &+ A_{1,\,n}\nabla q_{n}(\mathbf{x})\cdot \left( \nabla V_{n}(\mathbf{x})-\nabla \overline{q_{n-1}}(\mathbf{x})\right)\\
& = A_{2,n}|\nabla q_{n}(\mathbf{x})|^{2}+A_{3,\,n}\lvert\nabla V_n(\mathbf{x})-\nabla\overline{q_{n-1}}(\mathbf{x})\rvert^2,
\label{eq:q5}
\end{aligned}
\end{equation}
 where $A_{i,\,n}$, $i=1,\,2,\,3$, are some coefficients defined in \cite{BK}
 and can be computed analytically and  $\overline{q_{n-1}}%
 =h\sum_{j=0}^{n-1} q_{j}$.  The tail function $V=V_n$ is
 approximated iteratively, see algorithm below.  The discretized
 version of the boundary condition (\ref{eq:q2}) is given by
\begin{equation}
q_{n}(\mathbf{x})=\psi _{n}(\mathbf{x}):=\frac{1}{h}\int%
\limits_{s_{n}}^{s_{n-1}}\psi (\mathbf{x},\,s)\,\mathrm{d}s\approx \frac{1}{2}[\psi (%
\mathbf{x},\,s_{n})+\psi (\mathbf{x},\,s_{n-1})],\quad\mathbf{x}\in \partial \Omega .
\label{eq:q7}
\end{equation}%
 We also note that the first term on the right hand side of
 (\ref{eq:q5}) is negligible compared to the other terms since
 $\lvert A_{2,\,n}\rvert\sim \Lambda^{-1}$ for sufficiently
 large $\Lambda$, while $\lvert A_{i,\,n}\rvert \sim \Lambda^0$, $i=1,\,3$. Thus, we set
 $A_{2,\,n}|\nabla q_{n}|^{2}=0$. The system of elliptic equations
 (\ref{eq:q5}) with boundary conditions (\ref{eq:q7}) is solved
 sequentially starting from $n=1$. To solve it we use following
 algorithm:

\subsubsection*{Globally convergent algorithm}

\begin{itemize}
\item Compute the first tail  function $V_0$ (see \cite{BTKF} for details). Set
  $q_0 \equiv 0$.

\item For $n = 1,\, 2,\, \ldots,\, N$

\begin{enumerate}
\item Set $q_{n,\,0} = q_{n-1}$, $V_{n,\,1} = V_{n-1}$

\item For $i = 1,\, 2,\,\ldots,\, m_n$

\begin{itemize}
\item Find $q_{n,\,i}$ by solving (\ref{eq:q5})--(\ref{eq:q7}) with $%
V_{n}:=V_{n, \,i}$.

\item Compute $v_{n,\,i} = -h q_{n,\,i} - \overline{q_{n-1}} + V_{n,\,i}$.

\item Compute ${\varepsilon_{\mathrm{r},\,n,\,i}}$ via (\ref{eq:c}). Then solve the forward
problem (\ref{eq:fp1})--(\ref{eq:fp2}) with the new computed coefficient $%
{\varepsilon_\mathrm{r}} :={\varepsilon_{\mathrm{r},\,n,\,i}}$, compute $w:=w_{n,\,i}$
and update the tail $V_{n,\,i+1}$ by (\ref{eq:tail}).
\end{itemize}

\item Set $q_{n}=q_{n,\,m_{n}}$, ${\varepsilon_{\mathrm{r},\,{n}}}={\varepsilon_{\mathrm{r},\,n,\,m_{n}}}$, $V_{n}=V_{n,\,m_{n}+1}$ and go to the next frequency interval $\left[
s_{n+1},\,s_{n}\right] $ if $n<N.$ If $n=N$, then stop.
\end{enumerate}
\end{itemize}

Stopping criteria of this algorithm with respect to $i$ and $n$ are derived
computationally and is presented in \cite{BTKF,NBKF}. We denote the solution obtained at this stage by $\varepsilon_{\mathrm{r},\,\mathrm{glob}}$.

\section{Statement of Forward and Inverse Problems on the second stage}

\label{sec:stage2}

On the second stage we model the electromagnetic wave
propagation in an isotropic and non-magnetic space with permeability $\mu=1$ in
$\mathbb{R}^{3}$ with the dimensionless coefficient $\varepsilon
_\mathrm{r}$, which describes the spatially distributed
dielectric permittivity of the medium. We consider the following Cauchy
problem in the model problem for the electric field $E(\mathbf{x},\,t)= (E_1,\,E_2,\,E_3)(\mathbf{x},\,t)$
\begin{equation}\label{E_gauge1}
\begin{aligned}
&\varepsilon_\mathrm{r}(\mathbf{x}) \frac{\partial^2 E}{\partial t^2}(\mathbf{x},\,t) + \nabla \times
\big(\nabla \times E(\mathbf{x},\,t)\big) = (0,\,\delta(z-z_0)f(t),\,0),&&(\mathbf{x},\,t)\in\mathbb{R}^{3}
\times (0,\,T) , \\
&\nabla\cdot\big(\varepsilon_{\mathrm{r}}(\mathbf{x})E(\mathbf{x},\,t)\big) = 0, &&(\mathbf{x},\,t)\in \mathbb{R}^{3}\times (0,\,T) ,\\
&E(\mathbf{x},\,0) = 0, \quad \frac{\partial E}{\partial t}(\mathbf{x},\,0) =0, &&\mathbf{x}\in\mathbb{R}^{3}.
\end{aligned}
\end{equation}
 In the above equation $f\left( t\right)\not\equiv 0$ is the time-dependent waveform
of the incident plane wave. This wave  propagates along the $z$-axis and is incident at the plane $\left\{
z=z_{0}\right\}$.

We assume that the coefficient $\varepsilon _\mathrm{r}$ of
equation (\ref{E_gauge1}) is the same as in (\ref{2.20}).
 Let  again
$\Gamma \subset \partial \Omega $ be a part of the boundary $\partial
\Omega .$

\textbf{Coefficient Inverse Problem (CIP).} \emph{Suppose that the
  coefficient }$\varepsilon _\mathrm{r}
$\emph{ satisfies \eqref{E_gauge1}.  Determine the function
}$\varepsilon _\mathrm{r}\left( \mathbf{x}\right) $,
$\mathbf{x}\in \Omega $\emph{, assuming that the following function
}$g$\emph{
 is known for a single incident plane wave:}
\begin{equation}
E\left( \mathbf{x},\,t\right) =g\left( \mathbf{x},\,t\right) ~\forall \left(
\mathbf{x},\,t\right) \in \Gamma \times \left( 0,\,T\right) .  \label{2.5b}
\end{equation}

In (\ref{2.5b}) the function $g$ models time dependent
measurements of the electromagnetic field at the part $\Gamma $ of the
boundary $\partial \Omega $ of the domain $\Omega$ in which coefficient
$\varepsilon_\mathrm{r}$ is unknown. The uniqueness
of the above  CIP in the multidimensional case is currently known only if we
will consider in (\ref{E_gauge1}) a Gaussian function $\delta_{\theta}
\left(z-z_0\right) $ centered around $z_{0}$, which approximates the
function $\delta \left(z - z_0\right)$, or if at least one of initial
conditions in (\ref{E_gauge1}) is not zero. We again assume that uniqueness holds for our
CIP.

The function $E$ in (\ref{E_gauge1})
represents the voltage of one component of the electric field $E\left(
\mathbf{x},\,t\right) =\left( E_1,\,E_2,\,E_3\right) \left(
\mathbf{x},\,t\right)$. In our computer simulations of section \ref{sec:numex2}
the incident field has only one non-zero component $E_2$. This
component propagates along the $z$-axis until it reaches the target,
where it is scattered.
When solving the forward problem in our computations of section \ref{sec:numex2},
we first generate the data (\ref{2.5b}) by solving the problem
(\ref{E_gauge1}) for the case when the function
$\varepsilon_\mathrm{r}$ is taken as the one reconstructed by the
globally convergent method. Next, the computed component
$E_2$ on the surface $\Gamma$ is replaced with the measured
data. The other two components, $E_1$ and $E_3$,
are left the same as the ones obtained by the solution of the problem
(\ref{E_gauge1}), see details in \cite{BTKJ}.


\subsection{Domain decomposition
finite element/finite difference method}

To solve the problem (\ref{E_gauge1}) numerically we choose a bounded domain
$G$ such that $\Omega \subset G$.
In our computations of the second stage we use the domain
decomposition finite element/finite difference method of \cite{BM}.
To do that we decompose $G$ as $G =\Omega _\mathrm{FEM}\cup \Omega _\mathrm{FDM}$
with $\Omega _\mathrm{FEM}=\Omega$. Then, in computations, in $\Omega _\mathrm{FEM}$
a finite element method is used while in $\Omega _\mathrm{FDM}$ a finite
difference method is used, see details in \cite{BM}.

Using (\ref{2.20})  we have that
\begin{equation}
\begin{split}
\varepsilon _\mathrm{r}(\mathbf{x})& \geq 1,\text{ for }\mathbf{x}\in \Omega _\mathrm{FEM},
\\
\varepsilon _\mathrm{r}(\mathbf{x})& =1,\text{ for }\mathbf{x}\in \Omega _\mathrm{FDM}.
\end{split}
\notag
\end{equation}

As in \cite{BM} in our computations we used the following stabilized model problem
 with the parameter $\xi\geq 1$:
\begin{align}
&\varepsilon_\mathrm{r}(\mathbf{x}) \frac{\partial^2 E}{\partial t^2}(\mathbf{x},\,t) + \nabla \times \big( \nabla
\times E(\mathbf{x},\,t)\big)\notag\\ &\qquad\qquad\qquad- \xi\nabla \Big( \nabla \cdot\big(\varepsilon_\mathrm{r}(\mathbf{x}) E(\mathbf{x},\,t)\big)\Big) = 0,&&(\mathbf{x},\,t)\in G\times (0,\,T),  \label{model1_1} \\
&E(\mathbf{x},\,0) = 0,\quad \frac{\partial E}{\partial t}(\mathbf{x},\,0) = 0,&&\mathbf{x}\in G.
\label{model1_3}
\end{align}

To determine boundary conditions for \eqref{model1_1}, \eqref{model1_3}, we choose the domains $\Omega$ and $G$ such that
\begin{equation}
\Omega =\Omega _\mathrm{FEM}=\left\{ \mathbf{x}=\left(
x,\,y,\,z\right):-a<x<a,\,-b<y<b,\,-c<z<c'\right\} , \notag
\end{equation}
\begin{equation}
G=\left\{ \mathbf{x}=\left( x,\,y,\,z\right) :-A<x<A,\,-B<y<B,\,-C<z<z_{0}\right\} ,
\notag
\end{equation}
where $0 < a < A$, $0 < b < B$, $-C < -c < c' < z_0$,
and $\Omega _\mathrm{FDM}=G\setminus \Omega _\mathrm{FEM}.$ Denote by
\begin{equation}
\partial_{1}G :=\overline{G}\cap \left\{ z=z_{0}\right\} ,\quad
\partial_{2}G :=\overline{G}\cap \left\{ z=-C\right\} , \quad
\partial_{3} G:=\partial G\setminus \left( \partial _{1}G\cup \partial _{2}G\right) .
\notag
\end{equation}
The backscattering side of $\Omega$ is $\Gamma =\partial \Omega \cap \left\{
z=c'\right\}$. Next, define $\partial _{i}G_{T} :=\partial _{i}G\times
\left( 0,\,T\right)$, $i=1,\,2,\,3$. Let $t'\in \left( 0,\,T\right) $ be a number,
and we assume that the function $f\left( t\right) \in C\left[ 0,\,t'\right]$
and $f(t) = 0$ for $t > t'$.

Then boundary conditions for (\ref{model1_1})--(\ref{model1_3}) are:
\begin{align}
&E\left( \mathbf{x},\,t\right) = (0,\,f(t),\,0), &&(\mathbf{x},\,t)\in\partial_{1}G\times \left( 0,\,t'\right],  \label{70}\\
&\frac{\partial E}{\partial n}(\mathbf{x},\,t)=-\frac{\partial E}{\partial t}(\mathbf{x},\,t), &&(\mathbf{x},\,t)\in\partial _{1}G\times \left( t',\,T\right) ,  \label{7}\\
&\frac{\partial E}{\partial n} (\mathbf{x},\,t)=-\frac{\partial E}{\partial t}(\mathbf{x},\,t),&&(\mathbf{x},\,t)\in\partial _{2}G_{T},  \label{8}\\
&\frac{\partial E}{\partial n}(\mathbf{x},\,t)=0,&&(\mathbf{x},\,t)\in\partial _{3}G_{T},  \label{9}
\end{align}
where $\frac{\partial}{\partial n}$ is the normal derivative. Conditions (\ref{7}) and (\ref%
{8}) are first order absorbing boundary conditions \cite{EM}.
At the lateral boundaries we impose a homogeneous Neumann condition \eqref{9}.
In \cite{BM} it was shown that the solution to the original Maxwell's equations is well approximated by the solution to (\ref{model1_1})--(\ref{9}) in the case where $\xi=1$ and the discontinuities in $\varepsilon_\mathrm{r}$ are not too large.


The model problem (\ref{model1_1})--(\ref{9}) can be also rewritten as
\begin{align}
&\varepsilon_\mathrm{r}(\mathbf{x}) \frac{\partial^2 E}{\partial t^2}(\mathbf{x},\,t) + \nabla \big(\nabla \cdot E(\mathbf{x},\,t)\big) - \nabla \cdot \big(\nabla E(\mathbf{x},\,t)\big)\notag \\
&\qquad\qquad\qquad - \xi\nabla \Big( \nabla \cdot\big(\varepsilon_\mathrm{r}(\mathbf{x}) E(\mathbf{x},\,t)\big)\Big) = 0, &&(\mathbf{x},\,t) \in G\times (0,\,T),  \label{model3_1} \\
&E(\mathbf{x},\,0) = 0, \quad \frac{\partial E}{\partial t}(\mathbf{x},\,0)=0, &&\mathbf{x}\in G,\label{model3_2} \\
&E\left( \mathbf{x},\,t\right) = (0,\,f\left( t\right),\,0),&&(\mathbf{x},\,t)\in\partial _{1}G\times \left( 0,\,t'\right],  \label{model3_3} \\
&\frac{\partial E}{\partial n}(\mathbf{x},\,t) =-\frac{\partial E}{\partial t}(\mathbf{x},\,t),&&(\mathbf{x},\,t)\in\partial _{1}G\times \left( t',\,T\right),  \label{model3_4} \\
&\frac{\partial E}{\partial n} (\mathbf{x},\,t) = -\frac{\partial E}{\partial n}(\mathbf{x},\,t),&&(\mathbf{x},\,t)\in\partial _{2}G_{T},  \label{model3_5} \\
&\frac{\partial E}{\partial n}(\mathbf{x},\,t) = 0,&&(\mathbf{x},\,t)\in\partial _{3}G_{T}.\label{model3_6}
\end{align}
Here we have used the well-known  identity $\nabla \times (\nabla \times E) = \nabla (\nabla
\cdot E) - \nabla \cdot (\nabla E)$.
We refer to \cite{BM} for details of the numerical solution of the forward
problem (\ref{model3_1})--(\ref{model3_6}).

\subsection{Tikhonov functional}
\label{sec:invprobl}

We define
 $\Gamma'$
as the extension of the backscattering side $\Gamma $ up
to the boundary $\partial _{3}G$ of the domain $G$ that is,
\begin{equation}
\Gamma'=\left\{ \mathbf{x}=\left( x,\,y,\,z\right) : - X<x< X,\,- Y<y< Y,\,z=c'\right\} .\notag
\end{equation}
Let $G'$ be the part of the rectangular prism $G$ which lies between the
two planes $\Gamma'$ and $\{z = -C\}$:
\begin{equation}
G'=\left\{ \mathbf{x}=\left( x,\,y,\,z\right) :- X<x< X,\,- Y<y< Y,\,-C<z<c'\right\} .\notag
\end{equation}
Denote by $Q_{T}=G'\times \left( 0,\,T\right)$, and $S_{T}=\partial G'\times\left( 0,\,T\right)$.

In our CIP
 we have the data $g$ in (\ref{2.5b}) only on $\Gamma$.
These data are complemented on the rest of the
boundary $\partial G'$ of the domain $G'$ by simulated data using the immersing procedure of \cite{BTKJ}. Thus, we can  approximately get the function $\widetilde{g}$:
\begin{equation}
\widetilde{g}\left( \mathbf{x},\,t\right) = E\left( \mathbf{x},\,t\right),\quad\left( \mathbf{x},t\right) \in S_{T}.  \label{12}
\end{equation}

We solve  our inverse problem as an optimization problem.
To do so we minimize the Tikhonov functional:
\begin{equation}
F(E,\, \varepsilon_\mathrm{r}) := \frac{1}{2} \int_{S_T}\big(E(\mathbf{x},\,t) - \tilde{g}(\mathbf{x},\,t)\big)^2 z_{\delta}(t)\,\mathrm{d}\sigma\,\mathrm{d}t + \frac{1}{2} \gamma \int_{G}\big( \varepsilon_\mathrm{r}(\mathbf{x}) - {\varepsilon_{\mathrm{r},\,\mathrm{glob}}(\mathbf{x})}\big)^2\, \mathrm{d}\mathbf{x} ,  \label{functional}
\end{equation}
where $\gamma > 0$ is the regularization parameter and $\varepsilon
_{\mathrm{r},\,\mathrm{glob}} $ is the computed coefficient which
we have obtained on the first stage via the globally convergent
method.
Here, $z_{\delta}(t)$ is used to ensure the compatibility conditions at $%
\overline{Q}_{T}\cap \left\{ t=T\right\} $ for the adjoint problem, see \cite{BTKJ} for details of this function.

Let $E_\mathrm{glob}$ be the solution of the forward
problem (\ref{model3_1})--(\ref{model3_6}) with $\varepsilon _\mathrm{r}:=\varepsilon _{\mathrm{r},\,\mathrm{glob}}$. Denote
by $p=\frac{\partial E_\mathrm{glob}}{\partial n}\vert _{S_{T}}$. In addition to the Dirichlet condition (\ref{12}),
we set the Neumann boundary condition as
\begin{equation}
\frac{\partial E}{\partial n}\left( \mathbf{x},\,t\right) =p\left( \mathbf{x},\,t\right),
\quad \left( \mathbf{x},\,t\right) \in S_{T}. \notag
\end{equation}

Introduce the following spaces of real valued vector functions
\begin{equation*}
H_{E}^{1}(Q_{T})=\left\{ f\in [H^{1}(Q_{T})]^3:f(\mathbf{x},\,0)=0\right\} ,
\end{equation*}
\begin{equation*}
H_{\lambda }^{1}(Q_{T})=\left\{ f\in [H^{1}(Q_{T})]^3:f(\mathbf{x},\,T)=0\right\} ,
\end{equation*}
\begin{equation*}
U^{1}=H_{E}^{1}\left( G_{T}\right) \times H_{\lambda }^{1}\left(
G_{T}\right) \times B\left( G\right) ,
\end{equation*}%
where $B\left( G\right) $ is the space of functions bounded on $G$ with the
norm $\left\Vert f\right\Vert _{B\left( G\right) }=\sup_{G}\left\vert
f\right\vert .$ 

To minimize the functional (\ref{functional}) we introduce the Lagrangian
\begin{equation}  \label{lagrangian}
\begin{aligned}
L(E,\,\lambda,\, \varepsilon_\mathrm{r}) &= F(E,\, \varepsilon_\mathrm{r})
 - \int_{Q_T}\varepsilon_\mathrm{r}(\mathbf{x}) \frac{\partial \lambda }{\partial t}(\mathbf{x},\,t)\cdot \frac{\partial E}{\partial t}(\mathbf{x},\,t)\,\mathrm{d}\mathbf{x} \,\mathrm{d}t\\
&\quad - \int_{Q_T}\nabla \cdot E(\mathbf{x},\,t)\nabla \cdot\lambda(\mathbf{x},\,t)\,\mathrm{d}\mathbf{x}\,\mathrm{d}t
  + \int_{Q_T} \nabla E(\mathbf{x},\,t) \nabla \lambda(\mathbf{x},\,t)\,\mathrm{d}\mathbf{x}\,\mathrm{d}t\\
&\quad + \xi \int_{Q_T} \nabla \cdot \big(\varepsilon_\mathrm{r}(\mathbf{x}) E(\mathbf{x},\,t)\big)\nabla \cdot \lambda(\mathbf{x},\,t)\,\mathrm{d}\mathbf{x}\,\mathrm{d}t
 - \int_{S_T} \lambda(\mathbf{x},\,t)\cdot p(\mathbf{x},\,t)\,\mathrm{d}\sigma\,\mathrm{d}t,
\end{aligned}
\end{equation}
where $E$ and $\lambda$ are weak solutions of problems
(\ref{model4_1})--(\ref{model4_3}) and (\ref{16})--(\ref{18}),
respectively, see details in \cite{BTKJ}.

 We observe that in (\ref{lagrangian}) $ \left(E,\,\lambda,\,\varepsilon _\mathrm{r}\right)
 =w\in U^{1}$ and functions $E$ and $\lambda $ depend on the $\varepsilon _\mathrm{r}.$ To get the Fr\'{e}chet derivative
 $L^{\prime }$ of the Lagrangian (\ref{lagrangian}) rigorously, one
 should assume that variations of functions $E$ and $\lambda $ depend
 on variations of the coefficient $\varepsilon _\mathrm{r}$. It can be done
 similarly with section 4.8 of \cite{BK}. However for brevity here, to
 derive the Fr\'{e}chet derivative of the Lagrangian (\ref%
 {lagrangian}) we assume that in (\ref{lagrangian}) the elements of the vector
 function $(E,\,\lambda ,\,\varepsilon_\mathrm{r})$ can be varied independently
 of each other.

We search for a point $w \in U^1$ such that
\begin{equation}
L^{\prime }(w)\left( \overline{w}\right) =0,~~~\forall \overline{w}\in U^{1}.
\label{scalar_lagr}
\end{equation}
To find the Fr\'{e}chet derivative $L^{\prime }(w),$ we consider
$L\left( w+ \overline{w}\right) -L\left( w\right)$, for every
$\overline{w}\in U^{1}$ and single out the linear part, with respect to
$\overline{w}$, of the obtained expression.  Then the state
problem in the domain $G'$ is given by
\begin{align}
&\varepsilon_\mathrm{r}(\mathbf{x}) \frac{\partial^2 E}{\partial t^2}(\mathbf{x},\,t)
 + \nabla \big(\nabla \cdot E(\mathbf{x},\,t)\big) \notag\\&\qquad\qquad\qquad
 - \nabla \cdot \big(\nabla E(\mathbf{x},\,t)\big)
 - \xi\nabla \Big( \nabla \cdot\big(\varepsilon_\mathrm{r}(\mathbf{x}) E(\mathbf{x},\,t)\big)\Big) = 0, &&(\mathbf{x},\,t) \in Q_T,  \label{model4_1} \\
&E(\mathbf{x},\,0) = 0, \quad \frac{\partial E}{\partial t}(\mathbf{x},\,0) = 0,&& \mathbf{x}\in G', \label{model4_2} \\
&\frac{\partial E}{\partial n}\left( \mathbf{x},\,t\right) = p\left( \mathbf{x},t\right),&&(\mathbf{x},\,t)\in S_T.  \label{model4_3}
\end{align}
The adjoint problem is:
\begin{align}
&\varepsilon_\mathrm{r}(\mathbf{x}) \frac{\partial^2 \lambda}{\partial t^2}(\mathbf{x},\,t)
 + \nabla \big(\nabla \cdot\lambda(\mathbf{x},\,t)\big)\notag \\&\qquad\qquad\qquad
 - \nabla \cdot \big(\nabla \lambda(\mathbf{x},\,t)\big)
 - \xi\varepsilon_\mathrm{r}(\mathbf{x}) \nabla \big( \nabla \cdot \lambda(\mathbf{x},\,t)\big) =0,&&(\mathbf{x},\,t)\in Q_{T},  \label{16} \\
&\lambda (\mathbf{x},\,T)=0,\quad\frac{\partial \lambda}{\partial t}(\mathbf{x},\,T) =0,&&\mathbf{x}\in G', \label{17} \\
&\frac{\partial\lambda}{\partial t} (\mathbf{x},\,t) =z_{\delta }\left( t\right) \big(\widetilde{g}(\mathbf{x},\,t)- E(\mathbf{x},\,t)\big) \left( \mathbf{x},t\right),&&(\mathbf{x},\,t)\in S_{T}.\label{18}
\end{align}


\section{Finite element discretization}

\label{sec:fem}

For the finite element discretization of $\Omega_T = \Omega \times
(0,T)$ we used stabilized finite element method of \cite{BM}.  To do
that we define a partition $K_{h}=\{K\}$ of $G'$ which consists of
tetrahedra.  Here $h$ is a mesh function defined as $h|_{K}=h_{K}$ --
the local diameter of the element $K$. Let $J_{\tau }=\left\{
J\right\} $ be a partition of the time interval $(0,\,T)$ into
subintervals $J=(t_{k-1},\,t_{k}]$ of uniform length $\tau
  =t_{k}-t_{k-1}$. We also assume the minimal angle condition on the
  $K_{h}$ \cite{Brenner}.

To solve the state problem
(\ref {model4_1})--(\ref{model4_3}) and the adjoint problem
(\ref{16})--(\ref{18} ) we define the finite element spaces, $W_{h}^{E}\subset
H_{E}^{1}\left( Q_{T}\right) $ and $W_{h}^{\lambda }\subset H_{\lambda
}^{1}\left( Q_{T}\right) $. First, we introduce the finite element
trial space $W_{h}^{E}$ for every component of the electric field $E$
defined by
\begin{equation}
W_{h}^{E}:=\{w\in H_{E}^{1}(Q_T):w|_{K\times J}\in P_{1}(K)\times
P_{1}(J),~\forall K\in K_{h},~\forall J\in J_{\tau }\},  \notag
\end{equation}%
where $P_{1}(K)$ and $P_{1}(J)$ denote the set of linear functions on $K$
and $J$, respectively. We also introduce the finite element test space $%
W_{h}^{\lambda }$ defined by
\begin{equation}
W_{h}^{\lambda }:=\{w\in H_{\lambda }^{1}(Q_T):w|_{K\times J}\in P_{1}(K)\times
P_{1}(J),~\forall K\in K_{h},~\forall J\in J_{\tau }\}.  \notag
\end{equation}%
Hence, the finite element spaces $W_{h}^{E}$ and $W_{h}^{\lambda }$ consist
of continuous piecewise linear functions in space and time.
To approximate the function $\varepsilon _\mathrm{r}$, we use the space
of piecewise constant functions $V_{h}\subset L_{2}\left( \Omega \right) $,
\begin{equation}
V_{h}:=\{u\in L_{2}(\Omega ):u|_{K}\in P_{0}(K),~\forall K\in K_{h}\},
\notag
\end{equation}%
where $P_{0}(K)$ is the set of constant functions on $K$.


Next, we set $U_{h}=W_{h}^{E}\times W_{h}^{\lambda }\times V_{h}$.
The finite element method for solving equation (\ref{scalar_lagr}) now
reads: \emph{Find }$u_{h}\in U_{h}$\emph{, such that}
\begin{equation}
L^{\prime }(u_{h})(\bar{u})=0, ~\forall \bar{u}\in U_{h}. \notag
\end{equation}

\section{General framework for a posteriori error estimation for CIPs}

\label{sec:general}

 Let $(E_{h},\, \lambda_h,\, \varepsilon_h) \in U_h$ be
  finite element approximations of functions $(E,\, \lambda,\,
  \varepsilon_\mathrm{r}) \in U^1$, see details in \cite{BM, BMaxwell2}.
In our recent works \cite{BJ1, BJ,BK} we derived a posteriori error
estimates for three kinds of errors:
\begin{itemize}
\item  The error $|L(u) - L(u_h)|$ in the  Lagrangian  with $u=(E,\, \lambda,\, \varepsilon_\mathrm{r})$, and $u_h =(E_h,\, \lambda_h,\, \varepsilon_h)$. This error for hyperbolic CIPs was presented in \cite{BJ1, BJ}.

\item The error $|F({\varepsilon_\mathrm{r}}) - F(\varepsilon_h)|$ in the Tikhonov functional. This error for hyperbolic CIPs was derived in \cite{BK}.

\item The error $|{\varepsilon_\mathrm{r}} - \varepsilon_h|$ in the
  regularized solution of this functional ${\varepsilon_\mathrm{r}}$. This
  error for hyperbolic CIPs was presented in \cite{BK}.

\end{itemize}

To derive errors in the Lagrangian or in the Tikhonov functional we
first note that
\begin{equation}\label{errors}
\begin{split}
L(u) -  L(u_h) &= L^{\prime }(u_h)(u- u_h) +
 R(u,\,u_h), \\
F({\varepsilon_\mathrm{r}}) -  F(\varepsilon_h) &= F^{\prime }(\varepsilon_h)({\varepsilon_\mathrm{r}}- \varepsilon_h) +
 R(\varepsilon_\mathrm{r},\,\varepsilon_h),
\end{split}
\end{equation}
where $ R(u,\,u_h)$,and $R(\varepsilon_\mathrm{r},\,\varepsilon_h)$ are the second order remainders terms.  We
assume that $\varepsilon_h$ is located in the small neighborhood of $\varepsilon_\mathrm{r}$. Thus, the terms $R(u,\,u_h)$, $R(\varepsilon_\mathrm{r},\,\varepsilon_h)$
are small and we can neglect them in (\ref{errors}).

We now use the Galerkin orthogonality  principle
\begin{equation*}
\begin{split}
L^{\prime }(u_h)(\bar{u}) &= 0~~ \forall \bar{u} \in U_h,\\
 F^{\prime }(\varepsilon_h)(b) &= 0~~ \forall b \in V_h,
\end{split}
\end{equation*}
together with the splitting
\begin{equation*}
\begin{split}
u - u_h = (u - u_h^\mathrm{I}) + (u_h^\mathrm{I} - u_h), \\
\varepsilon_\mathrm{r} - \varepsilon_h = (\varepsilon_\mathrm{r} - \varepsilon_h^\mathrm{I}) + (\varepsilon_h^\mathrm{I} - \varepsilon_h),
\end{split}
\end{equation*}
where  $u_h^\mathrm{I}  \in U_h$ is the interpolant of $u$,
  and  $\varepsilon_h^\mathrm{I} \in V_h$  is the interpolant of $\varepsilon_\mathrm{r}$,
and get the following  representation of errors  in the Lagrangian and in the Tikhonov functional, respectively:
\begin{equation}\label{errorfunc}
\begin{split}
L(u) - L(u_h) &\approx L^{\prime}(u_h)( u - u_h^\mathrm{I}), \\
F(\varepsilon_\mathrm{r}) - F(\varepsilon_h) &\approx F^{\prime}(\varepsilon_h)( \varepsilon_\mathrm{r} - \varepsilon_h^\mathrm{I}).
\end{split}
\end{equation}
In the a posteriori error estimates (\ref{errorfunc}) we have two types of ``factors'':
\begin{itemize}
\item $L^{\prime}(u_h)$  and $F^{\prime}(\varepsilon_h)$ represent \emph{residuals}, and
\item  $u - u_h^\mathrm{I}$  and $\varepsilon_\mathrm{r} - \varepsilon_h^\mathrm{I}$  represent \emph{weights}.
\end{itemize}

The residuals of (\ref{errorfunc}) can be computed by knowing the
 finite element approximations $(E_{h},\, \lambda_h,\, \varepsilon_h)$, but the weights must be further estimated.


Let $f\in H^1(\Omega)$ be approximated by its piecewise linear interpolant $f_h^\mathrm{I}$ and finite element approximation $f_h$ over a mesh $K_h$ of $\Omega$ as outlined in Section \ref{sec:fem}. Standard interpolation estimates (following from, for instance, \cite{EEJ}) then gives
\begin{equation}
\left\| f-f_h^\mathrm{I}\right\| _{L_{2}\left( \Omega \right) }\leq  C_\mathrm{I}  \left\|h~ \nabla f\right\| _{L_{2}\left( \Omega \right) }.
\label{2.6}
\end{equation}
where $ C_\mathrm{I} = C_\mathrm{I}\left( \Omega,\,h \right) $ is positive
constant depending only on the domain $\Omega$ and the mesh function $h=h(x)$, the latter defined as in Section \ref{sec:fem}.
In addition, we can estimate  right hand side in (\ref{2.6}), see \cite{EEJ}, via
\begin{equation}
|\nabla  f| \leq \frac{|[f_h]|}{ h_{K}}, \label{2.6b}
\end{equation}
where $[f_h]$ denotes the normal jump of the function  $f_h$ over the edges of the element $K$.

 Similarly with (\ref{2.6}), \eqref{2.6b}  we estimate $u - u_h^\mathrm{I}$
 in terms of derivatives of the function $u$ and the mesh
parameters $h$ and $\tau$ as
\begin{equation}\label{est1}
| u - u_h^\mathrm{I} |\leq C_\mathrm{I} \left( h^2  \left|\frac{[u_h]_\mathrm{s}}{h} \right|  + \tau^2 \left|\frac{[u_h]_\mathrm{t}}{\tau} \right| \right),
\end{equation}
where $[u_h]_\mathrm{s}$ is the maximum modulus of a jump in
the normal derivative of $u_h$ across a side of the element $K$,
$[u_h]_\mathrm{t}$ is the maximum modulus of the jump of the time derivative
of $u_h$ across a boundary node of the time interval $J$, see details in \cite{BJ1, B, BMaxwell2, BJ}.

We also estimate $\varepsilon_\mathrm{r} -
\varepsilon_h^\mathrm{I}$ in terms of derivatives of the function
$\varepsilon_\mathrm{r}$ and the mesh parameter $h$ as
\begin{equation}\label{est2}
| \varepsilon_\mathrm{r} - \varepsilon_h^\mathrm{I} | \leq C_\mathrm{I} h \left | \frac{[\varepsilon_h]}{h} \right|.
\end{equation}
Here, $[\varepsilon_h]$ is the jump of the function $\varepsilon_h$
over the element $K$.
Substituting  estimates (\ref{est1}) and (\ref{est2})  in the right
hand side of (\ref{errorfunc}) we can compute a posteriori errors in
the Lagrangian or in the Tikhonov functional in explicit way as
\begin{equation}\notag
\begin{split}
| L(u) - L(u_h) | &\approx C_\mathrm{I} ||L^{\prime}(u_h)||\cdot( h ||[u_h]_\mathrm{s}|| + \tau ||[u_h]_\mathrm{t}||), \\
| F(\varepsilon_\mathrm{r}) - F(\varepsilon_h) | &\approx C_\mathrm{I} || F^{\prime}(\varepsilon_h)|| \cdot ||[\varepsilon_h]||.
\end{split}
\end{equation}

Finally,
to derive an estimate for the error $\varepsilon_\mathrm{r} - \varepsilon_h$ in the regularized solution
$\varepsilon_\mathrm{r}$ we use  the convexity property of the Tikhonov functional
 together with the interpolation property (\ref{2.6}).
Below we formulate theorem of \cite{BK} for the case of a posteriori
error estimate in the reconstructed function $\varepsilon_\mathrm{r}$ for the problem
(\ref{eq:fp1})--(\ref{eq:fp2}).

 \textbf{Theorem} \cite{BK} \emph{ \small
  Let $\varepsilon_h \in V_h$ be a finite element approximation of the
 solution $\varepsilon_\mathrm{r} \in H^1(\Omega)$ on the
  finite element mesh $K_h$ with the mesh function $h$. Then there
  exists a constant $D$ such that
$\left\| F^{\prime }\left( \varepsilon_1\right) - F^{\prime }\left(\varepsilon_2\right) \right\| \leq D\left\| \varepsilon_1- \varepsilon_2\right\|$ for every $\varepsilon_1$, $\varepsilon_2$ satisfying \eqref{2.20}.
  Then the following
  a posteriori error estimate for the regularized solution  $\varepsilon_\mathrm{r} $ holds}
\begin{equation*}
|| \varepsilon_h - \varepsilon_\mathrm{r} ||_{L^2(\Omega)} \leq   \frac{D}{\alpha } C_I  || h \varepsilon_h ||_{L_2(\Omega)}.
\end{equation*}

 \emph{Remark 5.1.} The natural question linked with the adaptivity is:
 Can one rigorously guarantee that the mesh obtained after the
 minimization of the Tikhonov functional on sequentially refined
 meshes of finite elements results in an improvement of the accuracy?
 For the first time this question was answered positively in
 \cite{BKK}, also, see the book \cite{BK} and the survey \cite{BK_AA}.

\section{Mesh refinement recommendation and the adaptive algorithm}

\label{subsec:ad_alg}


 In our adaptive algorithm for the mesh refinement we have used
ideas of \cite{BJ} and the Theorem 5.1 and criterion of the Remark 5.1
of \cite{BMaxwell2}.  From this criterion follows that the finite
element mesh should be locally refined in such subdomain of $\Omega $
where the maximum norm of the Fr\'{e}chet derivative of the objective
functional is large.

Define
\begin{equation}
\begin{aligned}
L_{h}^{\prime,\,m}(\mathbf{x})&=
 - {\int_{0}}^{T}\frac{\partial \lambda _{h}^{m}}{\partial t}(\mathbf{x},\,t)\cdot\frac{\partial E_{h}^{m}}{\partial t}(\mathbf{x},\,t)\,\mathrm{d}t\\&\quad
 + \xi\int_{0}^{T}\nabla\cdot E_{h}^{m}(\mathbf{x},\,t)\nabla \cdot \lambda _{h}^{m}(\mathbf{x},\,t)\,\mathrm{d}t
 + \gamma ({\varepsilon _{h}}^{m}(\mathbf{x})-{\varepsilon}_{\mathrm{r},\,\mathrm{glob}}(\mathbf{x})),
\label{Bhm}
\end{aligned}
\end{equation}
where $m$ is the iteration index in the optimization procedure, and
$(E_{h}^m,\, \lambda_h^m,\, \varepsilon_h^m)$ are finite element approximations of the functions $(E,\, \lambda,\, \varepsilon_\mathrm{r})$, see details in \cite{BM, BMaxwell2}.
\\

\textbf{Adaptive algorithm}

\begin{itemize}

\item  Step 0. Choose an initial mesh $K_{h}$ in $\Omega$ and an initial
time partition $J_{0}$ of the time interval $\left( 0,\,T\right) .$ Start from
the initial guess $\varepsilon_{h}^{0}= \varepsilon_\mathrm{r,\, glob}$.  Compute
the approximations $\varepsilon_{h}^{m}$ as:

\item Step 1. Compute the approximate solutions $E_{h}^m$ and $\lambda _{h}^m$ of the state problem (\ref{model3_1})--(\ref%
{model3_4}) and the adjoint problem (\ref{16})--(\ref{18}) on $K_{h}$ and $%
J_{k}$, using coefficient $\varepsilon_h^m$, and compute the Fr\'echet derivative $L^{\prime,\,m}_h$ via (\ref{Bhm}%
).

\item Step 2. Update the coefficient on $K_{h}$ using the
conjugate gradient method:
\begin{equation*}
\varepsilon_h^{m+1}(\mathbf{x}) : = \varepsilon_h^{m}(\mathbf{x}) + \alpha d^m(\mathbf{x}),
\end{equation*}
where $\alpha > 0$ is a step-size in the conjugate gradient method, and
\begin{equation*}
\begin{split}
d^m(\mathbf{x})&= -L^{\prime,\,m}_h(\mathbf{x}) + \beta^m d^{m-1}(\mathbf{x}),
\end{split}%
\end{equation*}
with
\begin{equation*}
\begin{split}
\beta^m &= \frac{|| L^{\prime,\,m}_h ||_{L_2(\Omega)}^2}{|| L^{\prime,\,m-1}_h ||_{L_2(\Omega)}^2},
\end{split}%
\end{equation*}
and $d^0(\mathbf{x}) = -L^{\prime,\,0}_h(\mathbf{x})$.

\item Step 3. Stop updating the coefficient  and set $\varepsilon_h : =
\varepsilon_h^{m+1}$, $M:= m+1$, if either $||L_h^{\prime,\,m}||_{L_{2}(
\Omega)}\leq \theta$ or norms $||\varepsilon_h^{m} ||_{L_{2}(\Omega)}$
are stabilized. Here $\theta$ is a tolerance number. Otherwise, set $m:=m+1 $
and go to step~1.

\item Step 4. Compute $L^{\prime,\,M}_h$ via (\ref{Bhm}). Refine the mesh at
all grid points $\mathbf{x}$ where
\begin{equation}
|L^{\prime,\,M}_h\left( \mathbf{x}\right) | \geq \beta _{1}\max_{\mathbf{x}\in\overline{%
\Omega }} |L_{h}^{\prime,M}\left( \mathbf{x}\right)|. \notag
\end{equation}
Here the tolerance number $\beta _{1}\in \left( 0,\,1\right) $ is chosen by
the user.

\item Step 5. Construct a new mesh $K_{h}$ in $\Omega$ and a new
partition $J_{k}$ of the time interval $\left( 0,\,T\right) $. On $J_{k}$ the
new time step $\tau $ should be chosen in such a way that the CFL condition
is satisfied. Interpolate the initial approximation $\varepsilon _\mathrm{r,\,glob}$
from the previous mesh to the new mesh. Next, return to step 1 at $m=1$ and
perform all above steps on the new mesh.
 Stop mesh refinements if norms defined in step 3 either
increase or stabilize, compared with the previous mesh.
\end{itemize}

 In step 2 of this algorithm $\alpha$ can be computed by
a line search procedure, see, for example, \cite{Peron}.

\section{Numerical studies}

\label{sec:numex}

In this section we present results of reconstruction of buried objects
placed inside a sand box using the two-stage numerical procedure.
 To do that we use the approximate globally
convergent algorithm of section \ref{sec:stage1} on the first stage and
the adaptive algorithm of section \ref{subsec:ad_alg} on the second
stage.

To collect experimental data we have used the same configuration as for the
targets placed in the air, see \cite{BTKF, NBKF} for details.  The
only difference is that in this work we consider the objects placed
inside a box filled with dry sand. The relative dielectric constant of dry sand
is $\varepsilon_\mathrm{r} \left( \text{sand}\right) =4$. We used this information
to model the case of buried objects.  In our experiment we have used
different types of targets, including both metallic and nonmetallic
ones. We refer to the Table 5.1 of \cite{NBKF2} for the full
description of all data sets. In this paper we present reconstruction of
four targets listed in the Table \ref{tab:table1}.  We refer to
\cite{NBKF2} for details of the data acquisition process.

In our computational studies  we had the following goals:

\begin{itemize}

\item to reconstruct refractive indices of dielectric targets and appearing dielectric constants of metals, and

\item to  image the location of targets, and their sizes and shapes.

\end{itemize}

To work with metallic objects, it is convenient to treat them as
dielectrics with large dielectric constants, see \cite{KBKSNF} for
details. We call these \emph{appearing
  dielectric constants} and choose values for them in the interval
\begin{equation}
\varepsilon _\mathrm{r}\left( \text{metallic target}\right) \in \left( 10,\,25\right).  \label{2.51}
\end{equation}%

 Using (\ref{2.51}), we set in all our tests the
upper value of the function $\varepsilon_\mathrm{r}$ as $%
b=25,$ see (\ref{2.20}). Thus, we set lower and upper bounds for the
reconstructed function $\varepsilon_\mathrm{r}$ in $\Omega$ as
\begin{equation}
M_{\varepsilon _\mathrm{r}}=\{\varepsilon _\mathrm{r}(\mathbf{x}):\varepsilon _\mathrm{r}\left(
\mathbf{x}\right) \in \left[ 1,\, 25\right],~\mathbf{x}\in\Omega \}.  \label{4.30}
\end{equation}
 We ensure the upper bound in (\ref{4.30}) via truncating to 25 those
 values of $\varepsilon_\mathrm{r}$ which exceed
 this number. Similarly we deal with the lower bound  of (\ref{4.30}).

To compare our computational results with directly measured refractive
indices $n=\sqrt{\varepsilon _\mathrm{r}}$ of dielectric targets and effective
dielectric constants of metallic targets (see (\ref{2.51})), we consider
the maximal values of the computed functions $\varepsilon_\mathrm{r}$ obtained in both algorithms,  and define
\begin{equation}
\varepsilon _\mathrm{r}^\mathrm{comp}=\max_{\mathbf{x}\in\overline{\Omega }}\varepsilon
_\mathrm{r}\left( \mathbf{x}\right) ,\quad n^{\mathrm{comp}}=\sqrt{\varepsilon _\mathrm{r}^\mathrm{%
{comp}}}. \label{4.300}
\end{equation}

\emph{Remark 7.1.}
As the objects we reconstruct are buried in dry sand with relative dielectric constant 4, our computational results should be scaled by that factor in order to obtain correct apparent dielectric constants and refractive indices. In Tables \ref{tab:table2}--\ref{tab:table6}, we present such scaled results.

\subsection{Data preprocessing}

We point out that
there is a \textit{huge misfit}
between our experimental data and computationally simulated data.
There are
several causes of this misfit listed in Section 4.2 of \cite{NBKF}.
 Because of this misfit, the \emph{central procedure} required before
 applying of our two-stage numerical procedure is data preprocessing.
 This procedure is heuristic and cannot be rigorously justified. In
 this work we have used the same data preprocessing procedure consisting of
 several steps as was used in \cite{NBKF, NBKF2}.  The three main steps in this
 data preprocessing are:

\begin{enumerate}
\item Data propagation.

\item Extraction of the targets signal from the total signal, which is a mixture of the signal from the target and the signal from the sand. This extraction is applied to propagated data.

\item Data calibration: to scale the measured data to the same
  scaling as in our simulations. In the case of the globally
  convergent method, a calibrating object was used. In the case of the
  above described adaptive finite element method a different calibration was used,
  see for details \cite{BTKJ}.
\end{enumerate}

We have propagated the data to a plane, which we call as the
\textit{propagated plane} and is located closer to the targets.  This
means that we approximate the scattered wave on the propagated plane
using the measured scattered wave on the measurement plane.  The
distance between the measurement plane and the target was found using
first time of arrival of the backscattered signal.
  Data calibration
is used to scale the measured data by a certain factor obtained in our
simulations.  We call this factor the \emph{calibration
  factor}. The choice of this factor is based on the data of a known
target which we call the \textit{calibrating object}.
The procedure of the extraction of the signal of the target from the total signal is more complicated and we refer to  \cite{NBKF2} for its many details.

\subsection{Computational domains}

We  choose our computational domain $G$  as
\begin{equation}\label{G}
G=\left\{ \mathbf{x=}(x,\,y,\,z)\in (-0.56,\,0.56)\times (-0.56,\,0.56)\times
(-0.16,\,0.1)\right\}.
\end{equation}
The boundary of the domain $G$ is $\partial
G=\partial _{1}G\cup \partial _{2}G\cup \partial _{3}G.$ Here, $\partial
_{1}G$ and $\partial _{2}G$ are front and back sides of the
domain $G$ at $\{z=0.1\}$ and $\{z=-0.16\}$, respectively, and $\partial _{3}G$ is the union of left, right, top and bottom
sides of this domain.

The the  domain $G$ is split
into two subdomains $\Omega_\mathrm{FEM}= \Omega $ and $\Omega_\mathrm{FDM}$ so that
$G= \Omega_\mathrm{FEM}\cup \Omega_\mathrm{FDM}$ and inner domain is defined as
\begin{equation}
\Omega_\mathrm{FEM}=\Omega =\left\{ \mathbf{x=}(x,\,y,\,z)\in (-0.5,\,0.5)\times
(-0.5,\,0.5)\times (-0.1,\,0.04)\right\} .  \label{8.0}
\end{equation}
The experimental data $g$ for both algorithms
are given at the front side $\Gamma$ of the domain
$\Omega$ which is defined as
\begin{equation}
\Gamma =\{\mathbf{x}=(x,\,y,\,z)\in \partial \Omega :z=0.04\} \notag
\end{equation}

In some tests of the first stage we used the shrunken computational domain
$G$  defined as
\begin{equation} \notag
 G =\left\{ \mathbf{x=}(x,\,y,\,z)\in (-0.24,\,0.24)\times
(-0.24,\,0.24)\times (-0.16,\,0.1)\right\} ,
\end{equation}
as well as the shrunken
computational domain
$\Omega_\mathrm{FEM}$  defined as
\begin{equation}  \label{zoom}
\Omega_\mathrm{FEM}=\Omega =\left\{ \mathbf{x=}(x,\,y,\,z)\in (-0.2,\,0.2)\times
(-0.2,\,0.2)\times (-0.1,\,0.04)\right\}.
\end{equation}

\subsection{Description of experimental data sets}

To test performance of both stages we have applied first the
approximate globally convergent algorithm and then an adaptive finite
element method to reconstruct the targets presented in Table
\ref{tab:table1}.  This table describes the details of used data sets
together with the burial depths of the targets. After obtaining computational results the refractive indices of all dielectric targets were measured, and these measured refractive indices were compared to those predicted by the computations.


Some of the non-blind targets were used
for the calibrating  procedure.  The
blind targets were used to ensure that our two-stage  procedure works in
realistic blind data cases.

We note that the burial depths of the targets of Table
\ref{tab:table1} varied between 3 cm to 5 cm. Typically
burial depths of antipersonnel land mines do not exceed 10 cm.
 The measured data of
the sand box  (without buried objects) was used for the
calibration of all data for the four objects of Table
\ref{tab:table1}.

\subsection{Numerical examples of the first  stage}

\label{sec:numex1}

In Tables \ref{tab:table2} and \ref{tab:table3} we summarize
reconstruction results for all objects  of Table \ref{tab:table1}. Table
\ref{tab:table2} shows shows the reconstructed refractive indices for the
non-metallic targets. For these targets, the refractive index
$n=\sqrt{\varepsilon_\mathrm{r} (\text{target})}$. Here, $\varepsilon_\mathrm{r} (\text{target})$ was
  chosen as $\varepsilon_\mathrm{r} (\text{target})=\max_{\mathbf{x}\in \Omega
  }\varepsilon_\mathrm{r} (\mathbf{x})$. Table
  \ref{tab:table3} shows the burial depths and the effective dielectric
  constants of the metallic targets.
From Tables \ref{tab:table2} and \ref{tab:table3} we can see that the burial depth was
accurately estimated in most cases, with the errors not exceeding 1 cm.

The estimates of the refractive indices of non-metallic targets with
refractive indices larger than that of the sand (water and wet wood) are
quite accurate with the average error of about 8.5\%.

Note that the error in our direct measurement of the refractive index
of the wet wood was 10\%. For water, we were unable to directly
measure its refractive index at the used frequency of the signal,
which was about 7.5 GHz. Therefore, we have made a separate experiment
described in \cite{NBKF2} where we have obtained a reference value  $n=4.88$ for
water.  We observe from Table \ref{tab:table2} that for
water we have obtained a value of $n$ close to the reference value.
Targets with smaller refractive indices than that of the sand are
modelling plastic land mines and improvised explosive devices
(IEDs). We have observed that in this case we can image these targets
only if their burial depths do not exceed 5 cm, see for example,
reconstruction of target 3 in Table \ref{tab:table2} and in Figure
\ref{fig:1}-c).

In our experiments we observed that the signals of the metallic
targets were stronger compared to the signal from sand.  In our
previous works, we have established that the effective dielectric
constant of metals should be larger than 10--15, see \cite{BTKF,
  NBKF}.  From Table  \ref{tab:table3}  we see that we have obtained similar results as
in our previous studies.

From Table \ref{tab:table1} we observe that in our experiments we were
supposed to reconstruct two metallic blocks which were placed at 1 cm
separation to each other. On the other hand, the wavelength $\lambda$
of our device is 4.5 cm. Thus, $\lambda/4.5$ is the
distance between these two targets and \textit{superresolution}
is achieved beyond the diffraction limit.  Table \ref{tab:table3} and
Figure \ref{fig:1}-d) shows that we have accurately imaged both
targets. 
This phenomenon was not expected and should be studied
further because of its importance when combined with quantitative
imaging.

\subsection{Numerical examples of the second stage}

\label{sec:numex2}

  From the results of the first stage we can conclude that this stage
  provides accurate locations of the targets as well as accurate values of the
  refractive indices $n=\sqrt{\varepsilon _\mathrm{r}}$ of the dielectric targets
  and large values of effective dielectric constants $\varepsilon
  _\mathrm{r}$ for the metallic targets of interest. However, the approximate
  globally convergent algorithm does not reconstruct the
  shapes of the targets in the $z$-direction well, see Figure \ref{fig:1}. Because of that we have used the second stage
  where we have minimized  the Tikhonov functional on locally
  adaptively refined meshes.

\subsubsection{Computations of the forward problem}

The data $g$ in our
experiments of the second stage are given only for the second
component $E_2$  of the
electric field $E$ in (\ref{2.5b}) and are measured at the front side
$\Gamma$ of the domain $\Omega$ which is defined as
\begin{equation}
\Gamma =\{\mathbf{x}=(x,\,y,\,z)\in \partial \Omega :z=0.04\}.  \notag
\end{equation}

To generate backscattering data for other two components $E_1$ and $E_3$ we
solve the forward problem (\ref{model3_1})--(\ref{model3_6}) in the
computational domain $G$ defined as in the first stage in (\ref{G})
with the known value of $\varepsilon_\mathrm{r}$ obtained at the first stage of
our two-stage numerical procedure. We use a stabilized domain
decomposition method of \cite{BM} implemented in the software package
WavES \cite{waves}. We split $G$ into two subdomains $\Omega_\mathrm{FEM}= \Omega
$ and $\Omega_\mathrm{FDM}$ so that $G=\Omega_\mathrm{FEM}\cup \Omega_\mathrm{FDM}$ and the inner domain is
defined as in (\ref{8.0}).

 Once the forward problem (\ref{model3_1})--(\ref{model3_6}) is solved
  to generate backscattering data  for the two components $E_1$ and $E_3$ at the boundary
 $\Gamma'$, then after the data immersing procedure described in  Section 7.3.3 of \cite{BTKJ} the
 inverse problem is solved via the algorithm of section
 \ref{subsec:ad_alg}.
The immersing procedure of \cite{BTKJ} immerses the time-dependent
propagated experimental data $g\left( \mathbf{x} ,\,t\right)
=E_{2}\left( \mathbf{x},\,t\right) \rvert _{\mathbf{x}\in \Gamma}$ into the
computationally simulated data and then extends the data
$g$ from $\Gamma $ to $\Gamma'$.

We choose the waveform
 $f$ in (\ref{model3_1})--(\ref{model3_6}) as
\begin{equation*}
f(t)=\sin (\omega t),\quad0\leq t\leq t':=\frac{2\pi }{\omega },
\end{equation*}%
where we use $\omega =30$ and $T=1.2.$ We solve the problem (\ref{model3_1})--(\ref{model3_6}) using the
explicit scheme of \cite{BM}  with the time step size $\tau =0.003$, which satisfies the
CFL condition.

\subsubsection{Reconstructions}

\label{sec:8.2}

Suppose that in the adaptive  algorithm of section
\ref{subsec:ad_alg} we have obtained the function
$\varepsilon_\mathrm{r}$.  We obtain then the image of the
dielectric targets based on the function $\varepsilon_\mathrm{r,\,diel}$ which we define as
\begin{equation*}
\varepsilon _\mathrm{r,\,diel}\left( \mathbf{x}\right) = \left\{
\begin{array}{l}
\varepsilon_\mathrm{r}\left( \mathbf{x}\right) \text{ if } \varepsilon_\mathrm{r}\left( \mathbf{x}\right) \geq 0.5\max_{\mathbf{x}\in\overline{\Omega }}%
\varepsilon_\mathrm{r}\left( \mathbf{x}\right) , \\
1\text{ otherwise.}
\end{array}
\right.
\end{equation*}
For  metallic targets we used similar function
 $\varepsilon_\mathrm{r,\,metal}$,
\begin{equation*}
\varepsilon _\mathrm{r,\,metal}\left( \mathbf{x}\right) =\left\{
\begin{array}{l}
\varepsilon_\mathrm{r}\left( \mathbf{x}\right) \text{ if } \varepsilon_\mathrm{r}\left( \mathbf{x}\right) \geq 0.5 \max_{\mathbf{x}\in\overline{\Omega }}%
\varepsilon_\mathrm{r}\left( \mathbf{x}\right) , \\
1\text{ otherwise.}
\end{array}
\right.
\end{equation*}

In our experiments we apply the adaptive algorithm of section
\ref{subsec:ad_alg} to improve shape of targets listed in Table
\ref{tab:table1}.

\begin{figure}[tbp]
\begin{center}
\begin{tabular}{cc}
{\includegraphics[scale=0.45,clip=]{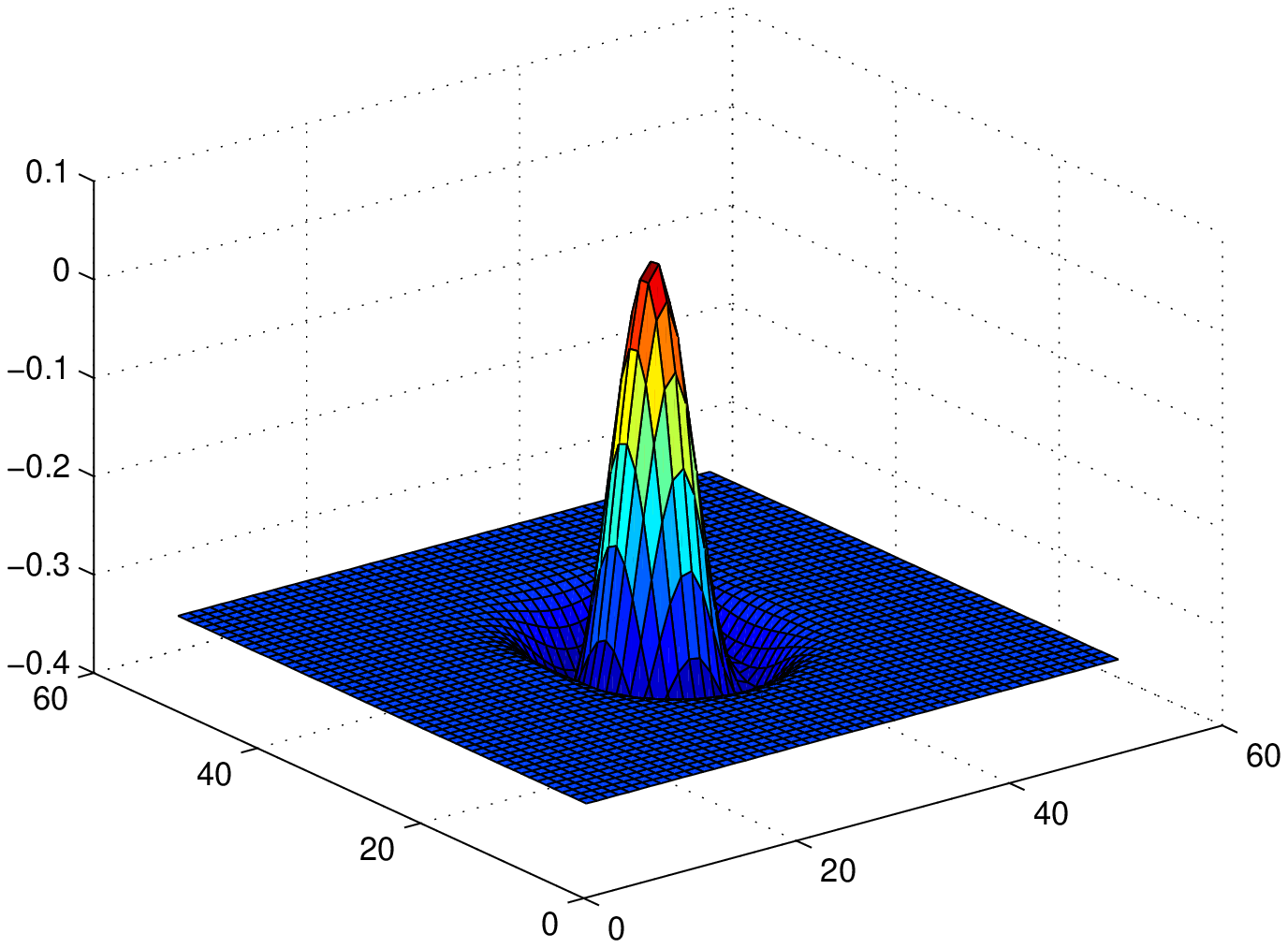}} &
{\includegraphics[scale=0.45, clip=]{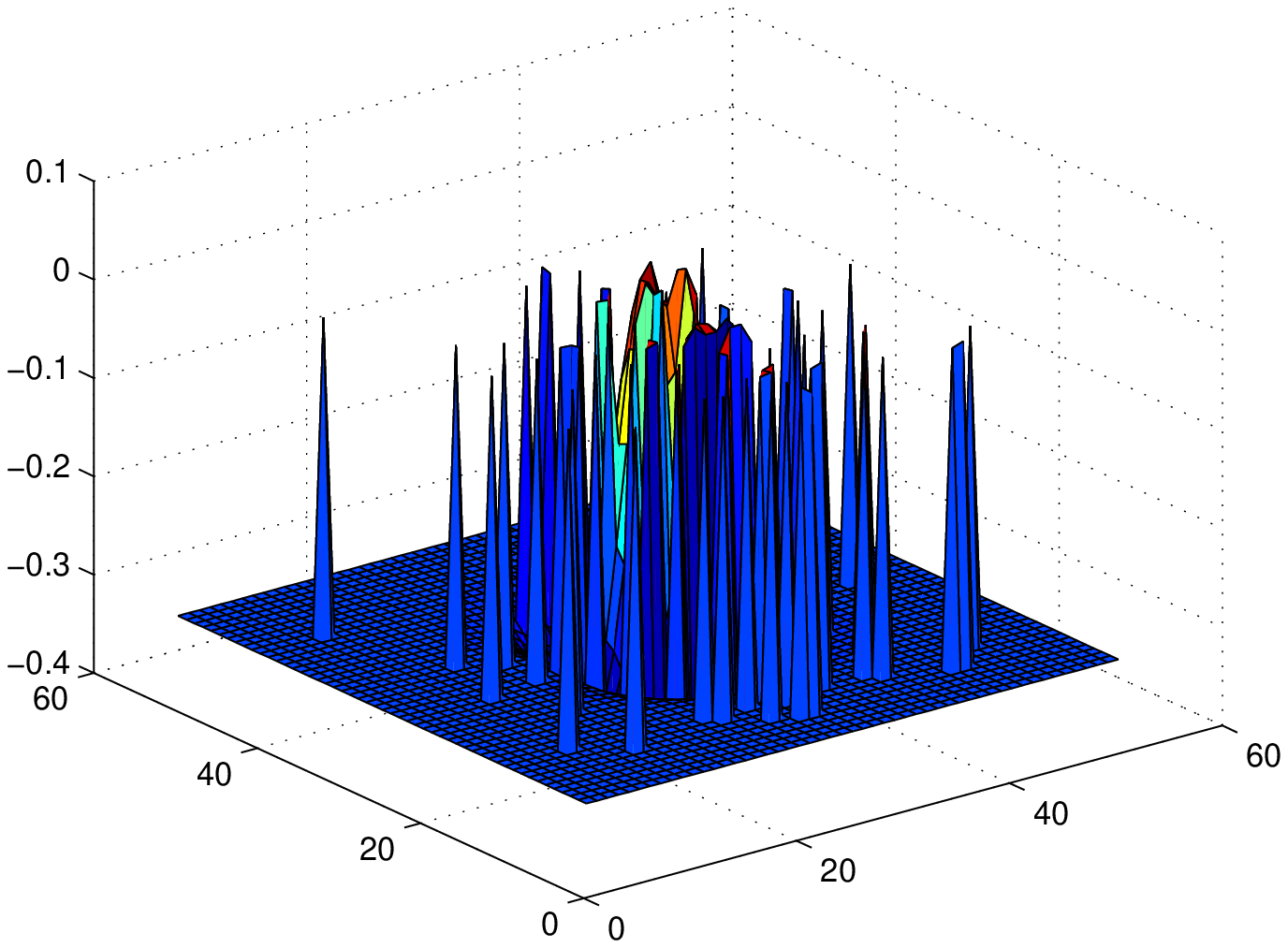}} \\
a) t=0.3 & b)   t=0.3 \\
{\includegraphics[scale=0.45, clip=]{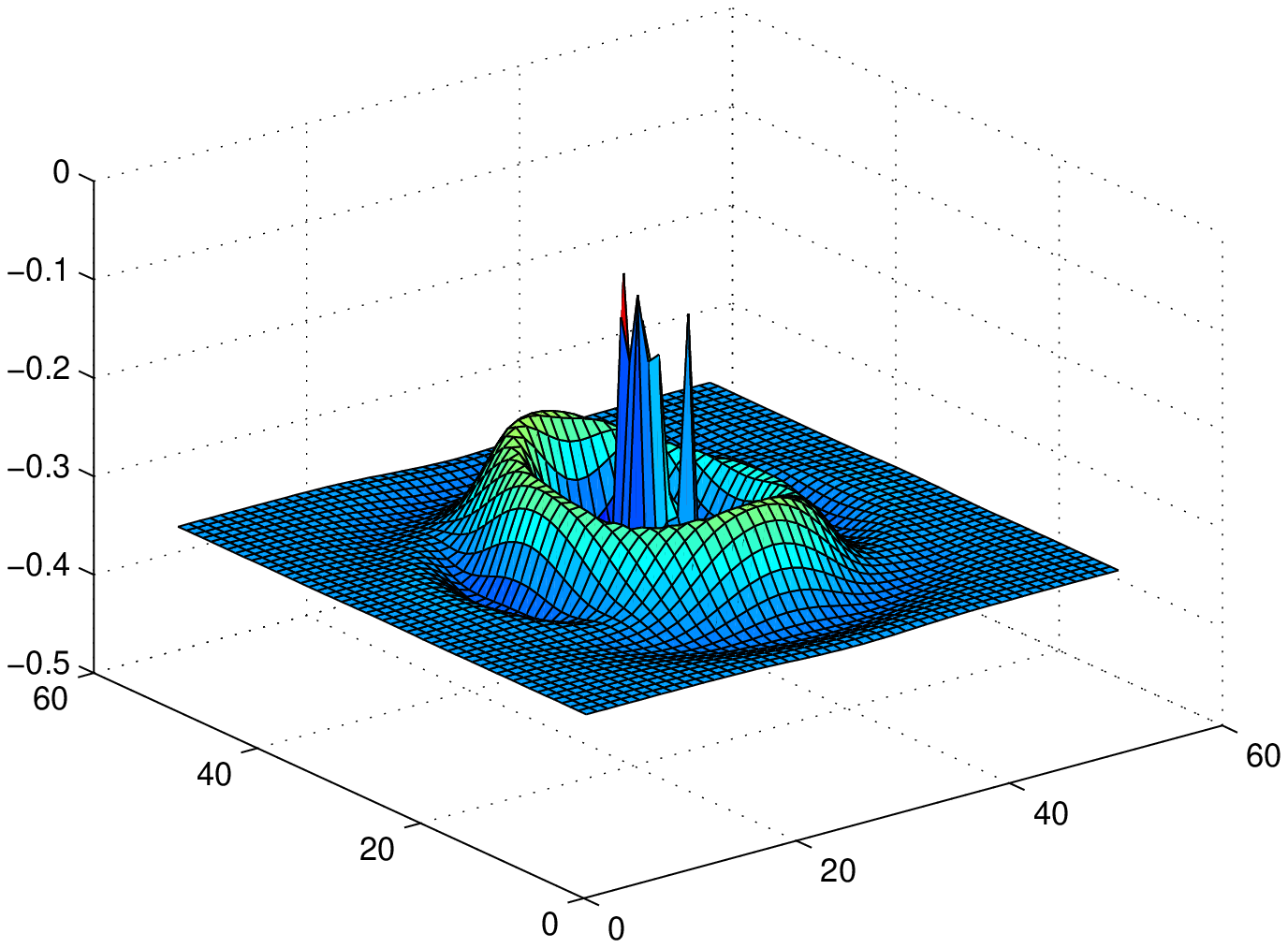}} &
{\includegraphics[scale=0.45, clip=]{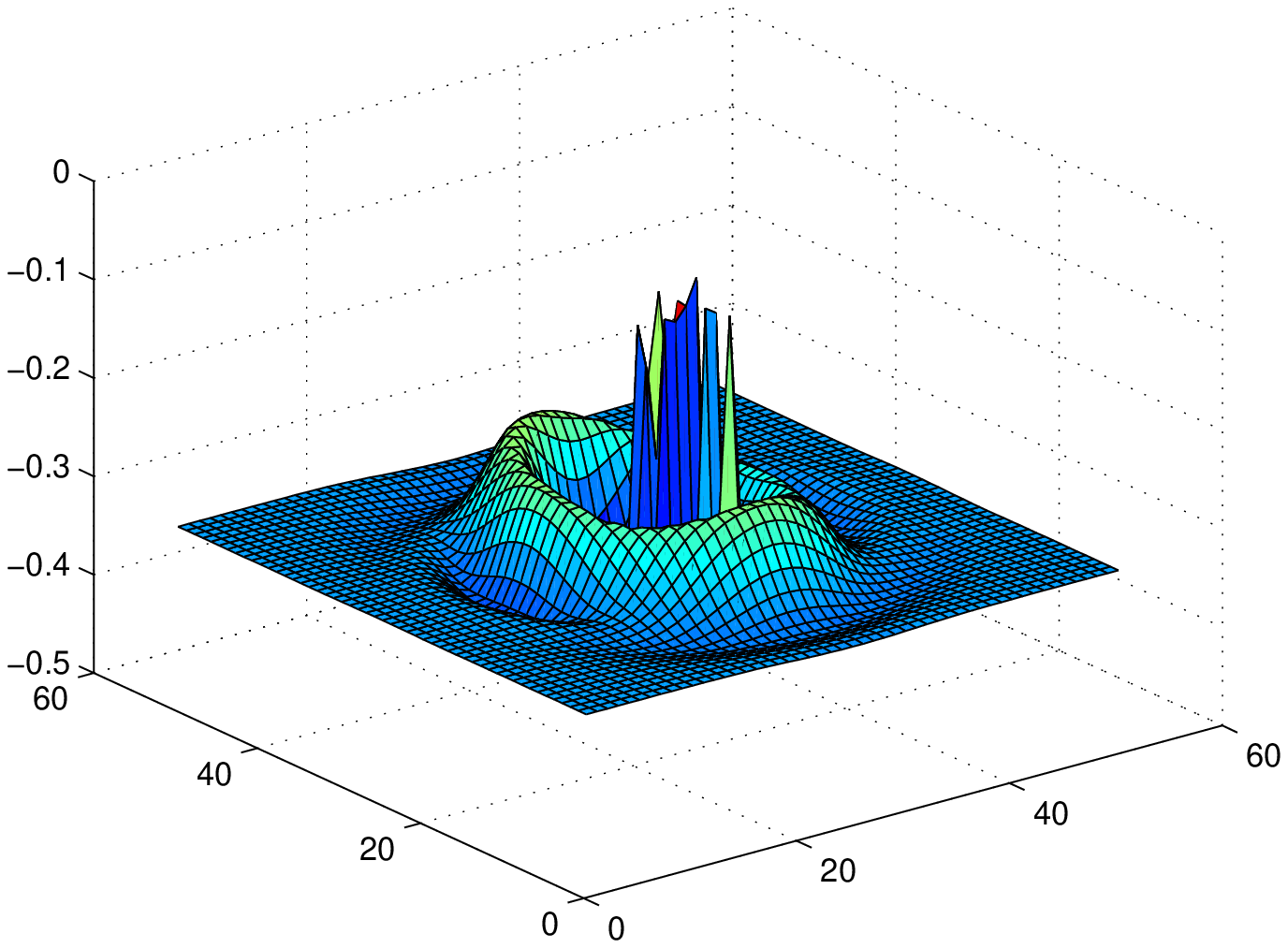}} \\
c)  t=0.45 & d) t=0.45
\end{tabular}
\end{center}
\caption{Backscattering immersed data of the second component $E_2$ of
the  electric field for object 4 (two metallic blocks at 1 cm separation)
  of Table \protect\ref{tab:table1}.
 On the left we
  show backscattering immersed data which are immersed into measured data without presence of sand,
  on the right - with presence of sand. Recall that the final time is
  $T=1.2$.}
\label{fig:4}
\end{figure}

\begin{figure}[tbp]
\begin{center}
\begin{tabular}{cc}
{\includegraphics[scale=0.25, angle=-90,clip=]{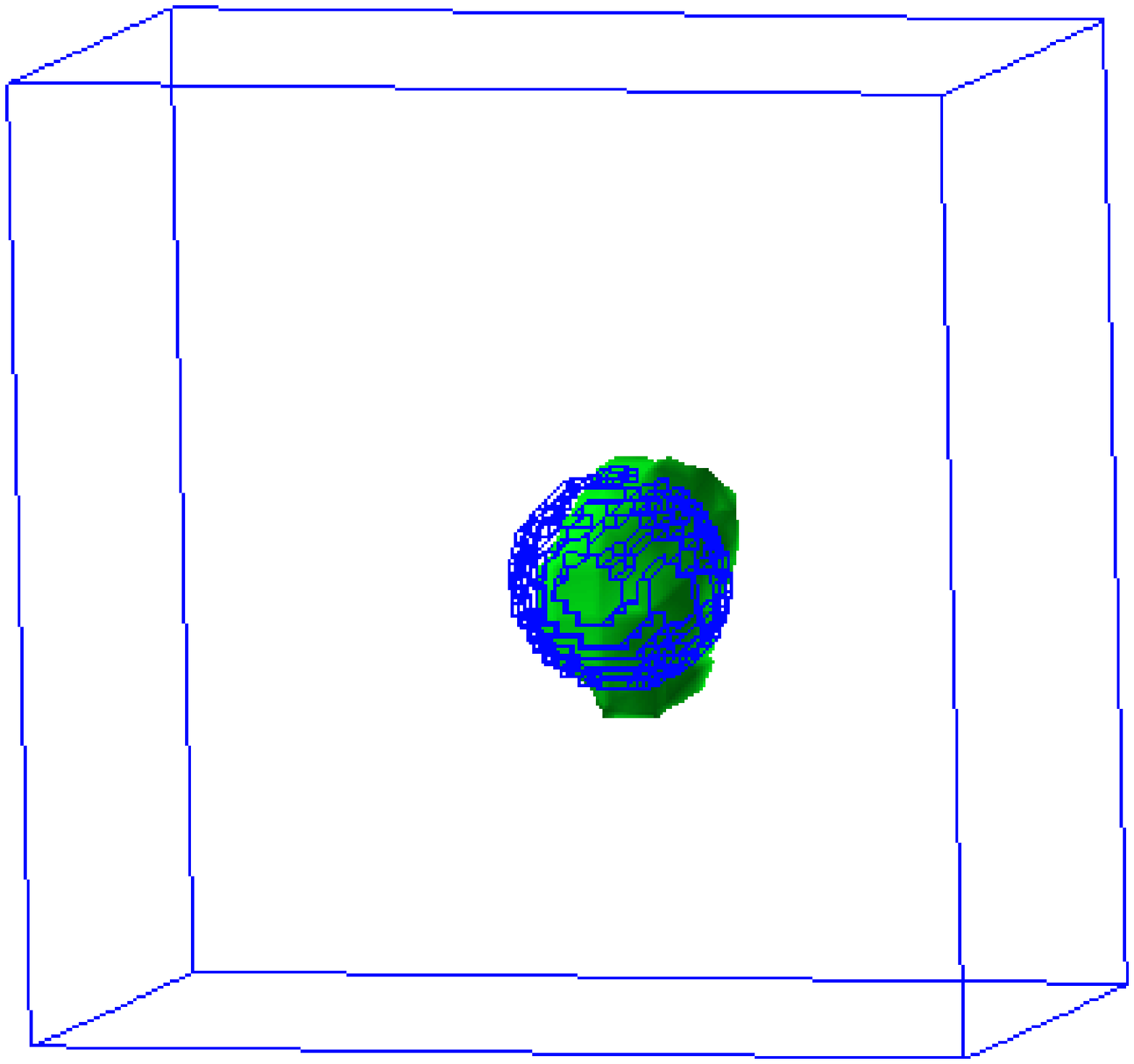}} &
{\includegraphics[scale=0.25,angle=-90, clip=]{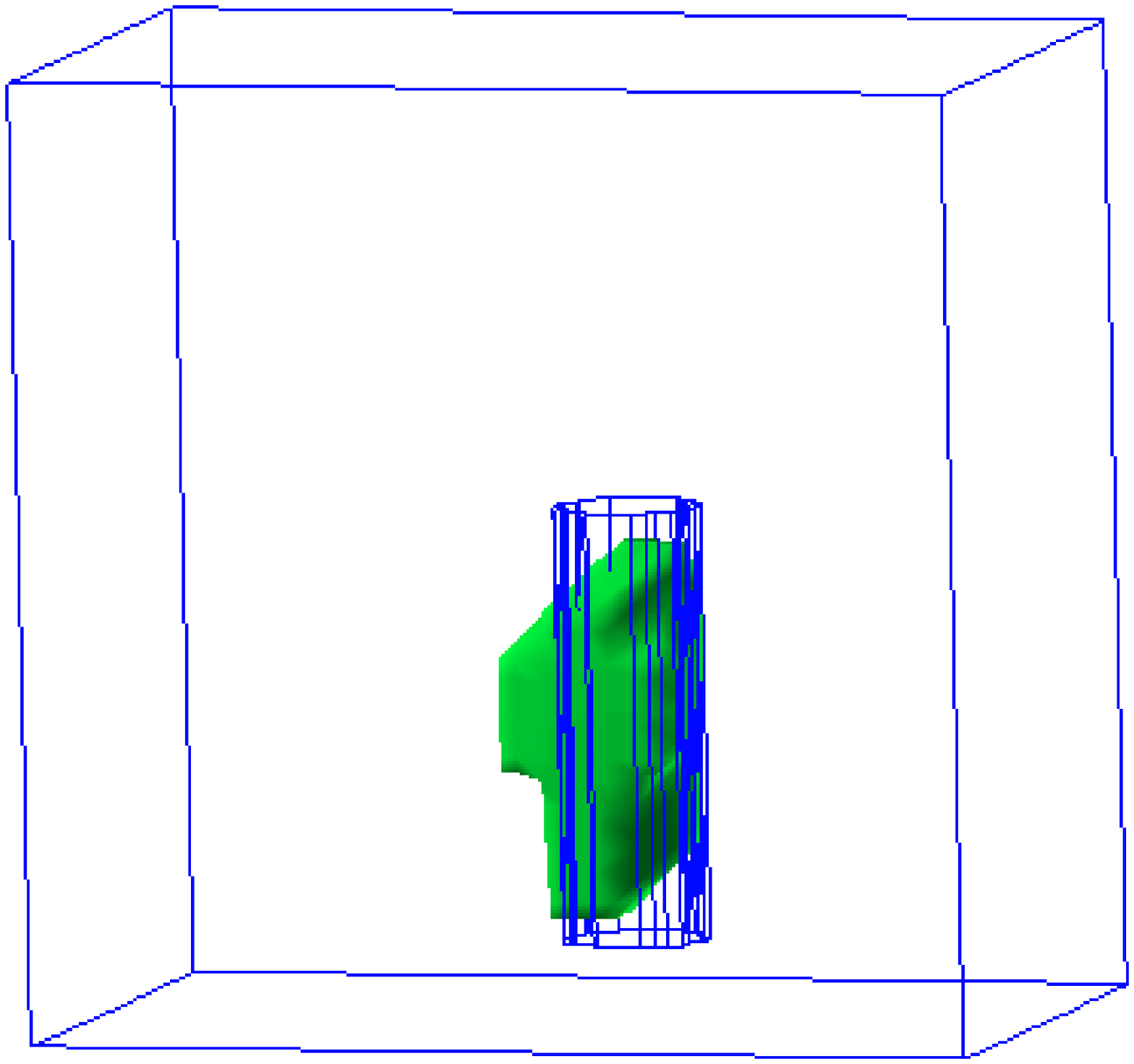}} \\
a)  target 1 & b)  target 2 \\
{\includegraphics[scale=0.25, angle=-90,clip=]{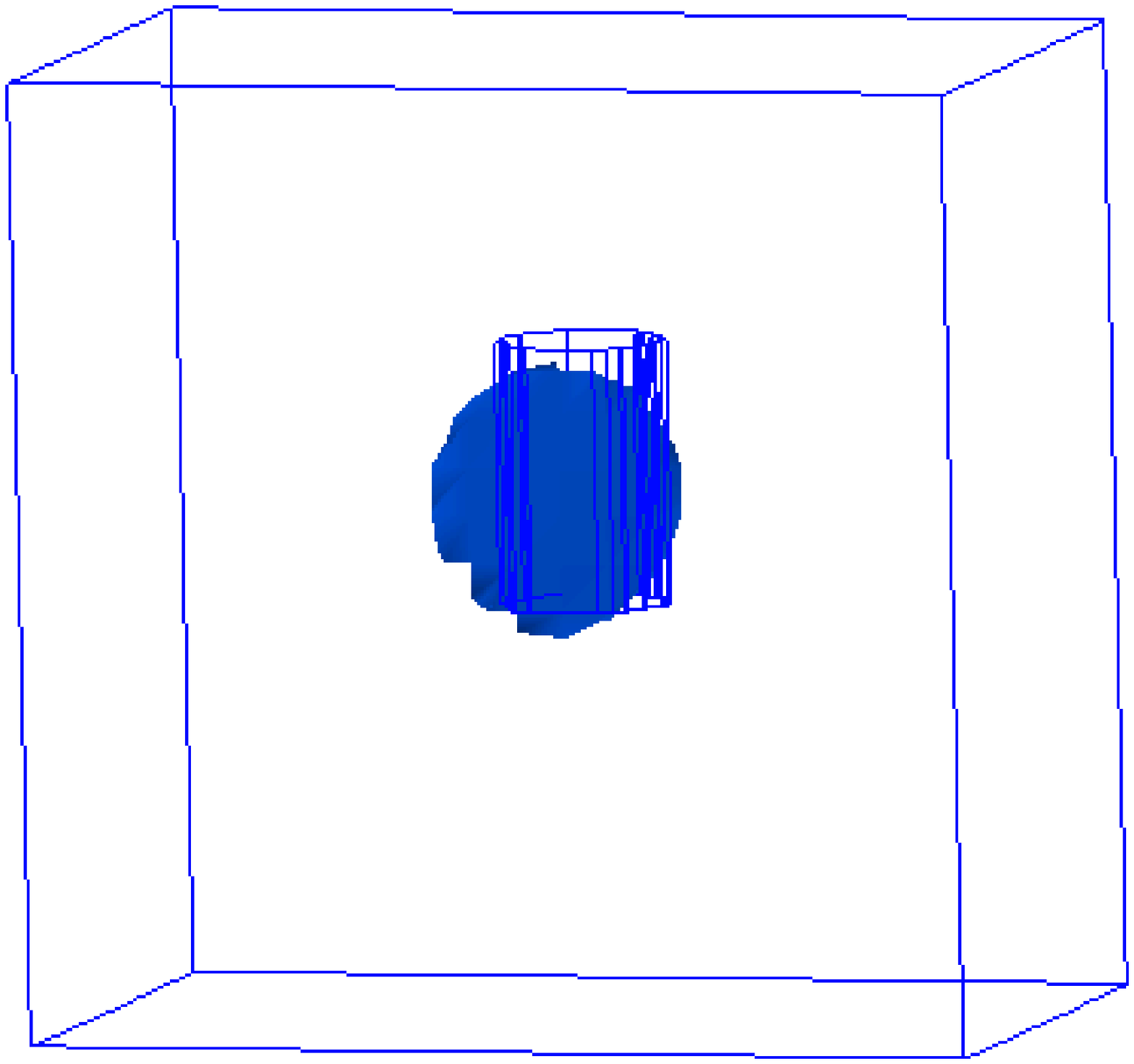}} &
{\includegraphics[scale=0.25,angle=-90, clip=]{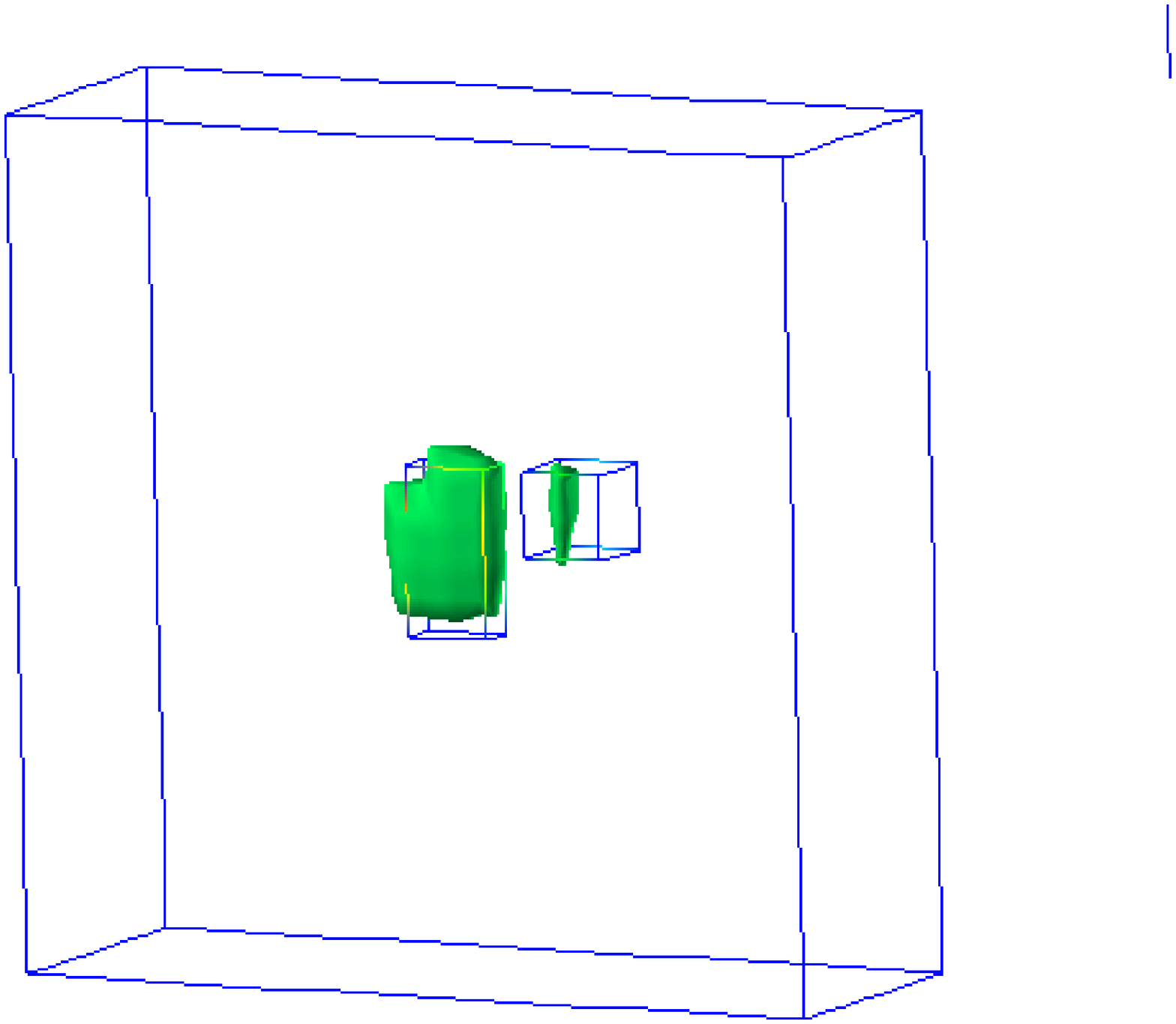}} \\
c) target 3 & d) targets 4
\end{tabular}
\end{center}
\caption{{\protect\small \emph{Reconstructions of  targets of Table \protect\ref{tab:table1}
obtained on the first stage of our two-stage
numerical procedure.}}}
\label{fig:1}
\end{figure}

Recall that to apply immersing procedure of the experimental data
$g$ into simulated data $E_2$ we solve
the problem  (\ref{model3_1})--(\ref{model3_6}) numerically with the known values of the function
$\varepsilon_\mathrm{r}=\varepsilon_\mathrm{r,\,glob}$ obtained at the first stage of our two-stage numerical
procedure, see Tables \ref{tab:table2}, \ref{tab:table3} for the
function $\varepsilon_\mathrm{r,\,glob}$.  Figure \ref{fig:4} show
backscattering immersed data of the second component of electric field
$E_2$ for target \#4 (two metallic blocks) of Table \ref{tab:table1} at
different times.

Table \ref{tab:table5} lists both computed refractive index $n^{\mathrm{comp}}$, obtained via (\ref{4.300}), on adaptively refined meshes and directly measured refractive indices $n$ of the dielectric targets. Table~\ref{tab:table6} lists calculated appearing dielectric constants
$\varepsilon _\mathrm{r}^{\mathrm{comp}}$ of the metallic targets.
From Table \ref{tab:table6} we observe that $\varepsilon _\mathrm{r}^{\mathrm{comp}}>10$ for all
metallic targets, and thus   (\ref{2.51}) is satisfied.

An important observation, which can be deduced from Table
\ref{tab:table6}, is that our adaptive algorithm can still compute large
inclusion/background contrasts exceeding 10:1.

\begin{figure}[tbp]
\begin{center}
\begin{tabular}{cc}
{\includegraphics[scale=0.25, angle=-90,clip=]{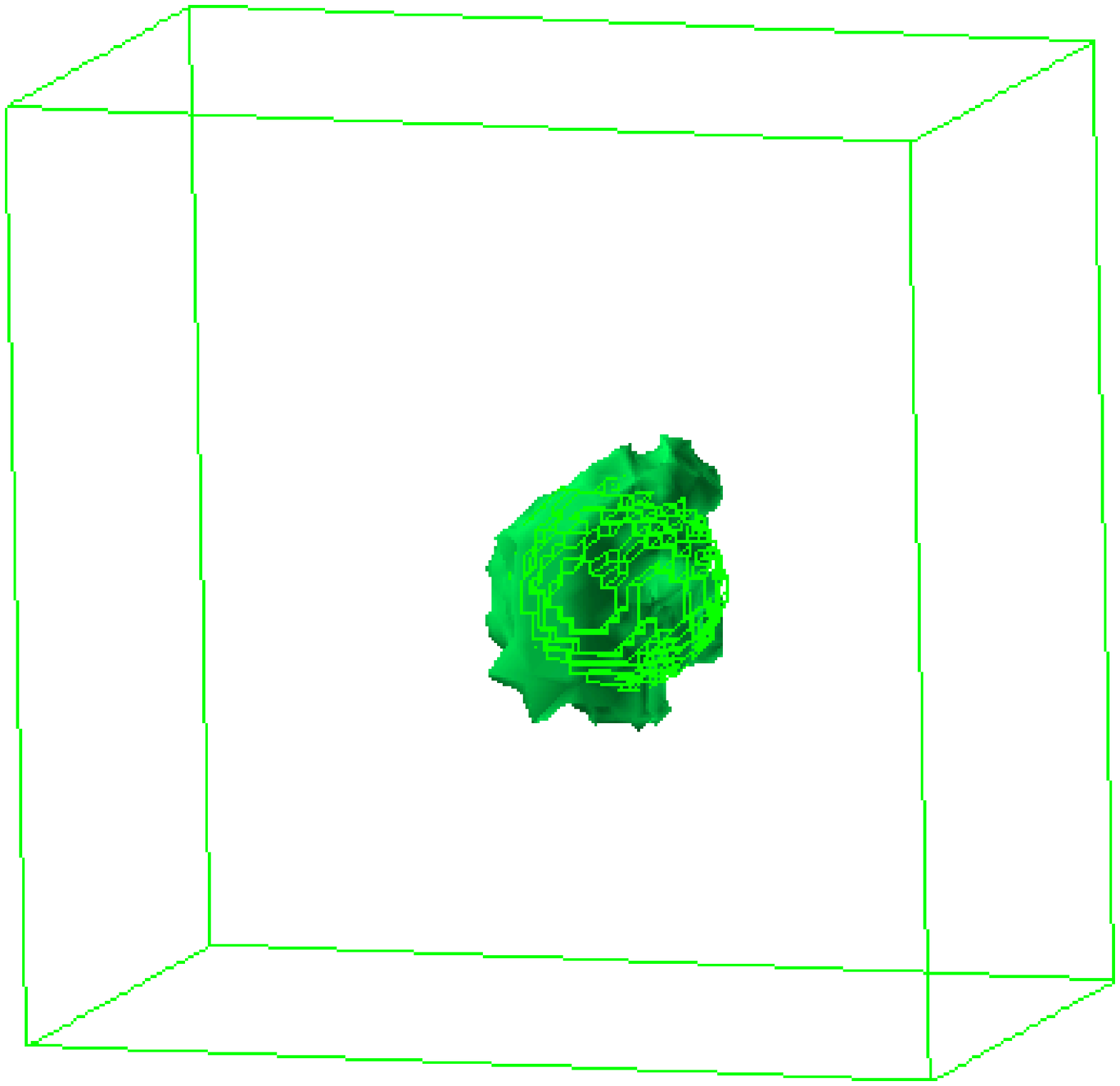}} &
{\includegraphics[scale=0.25,angle=-90, clip=]{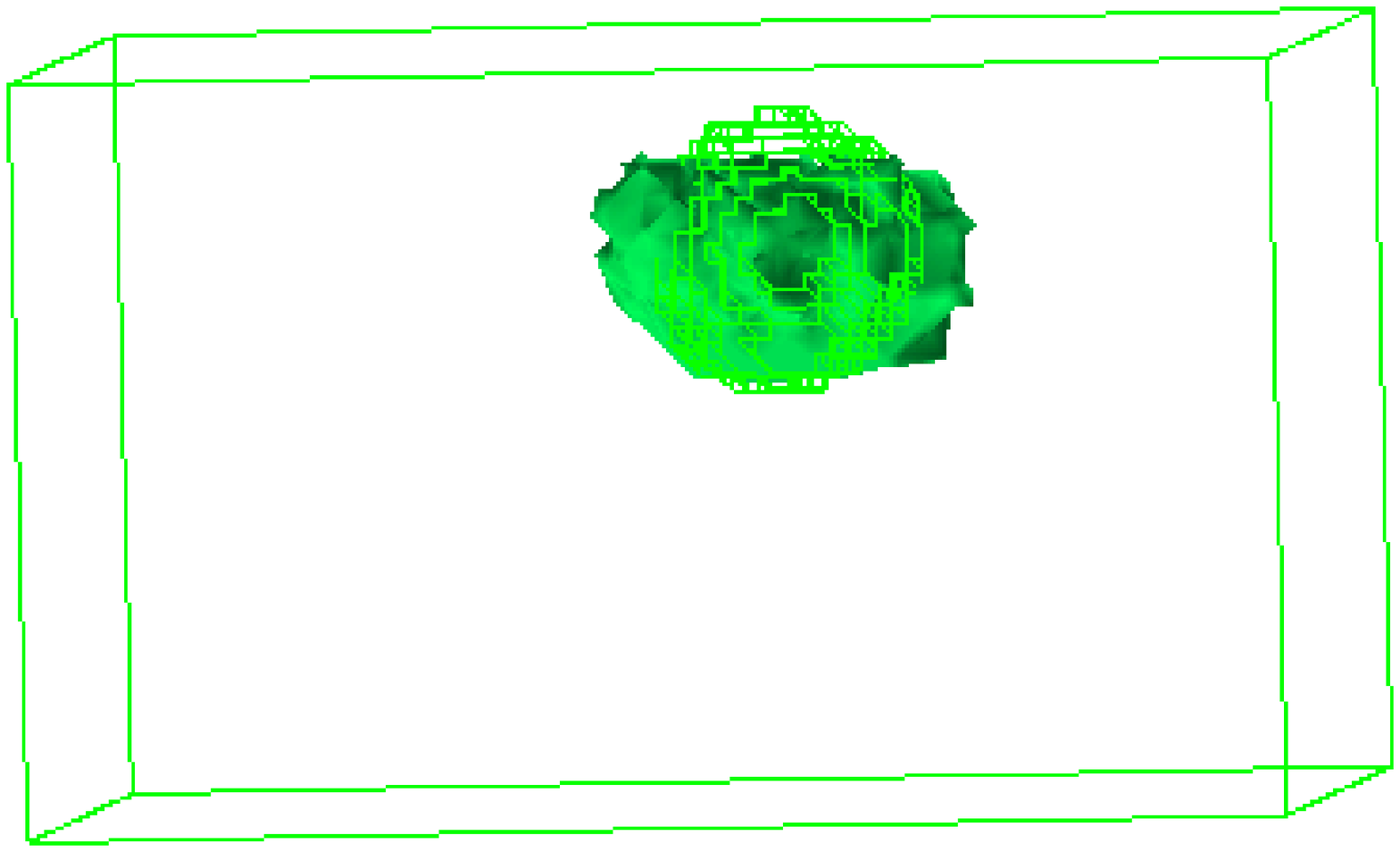}} \\
a) three times refined mesh, $xy$-view & b)  three times refined mesh, $yz$-view \\
\end{tabular}
\end{center}
\caption{{\protect\small \emph{
Computed image of target number 1 of Table \ref{tab:table1}
.
Thin lines indicate correct
shape. To have a better visualization we have zoomed the domain $\Omega$ in
(\protect\ref{8.0}) in the domain $\Omega_\mathrm{FEM}$ in (\protect\ref{zoom}).
 }}}
\label{fig:2}
\end{figure}

\begin{figure}[tbp]
\begin{center}
\begin{tabular}{cc}
{\includegraphics[scale=0.25, angle=-90,clip=]{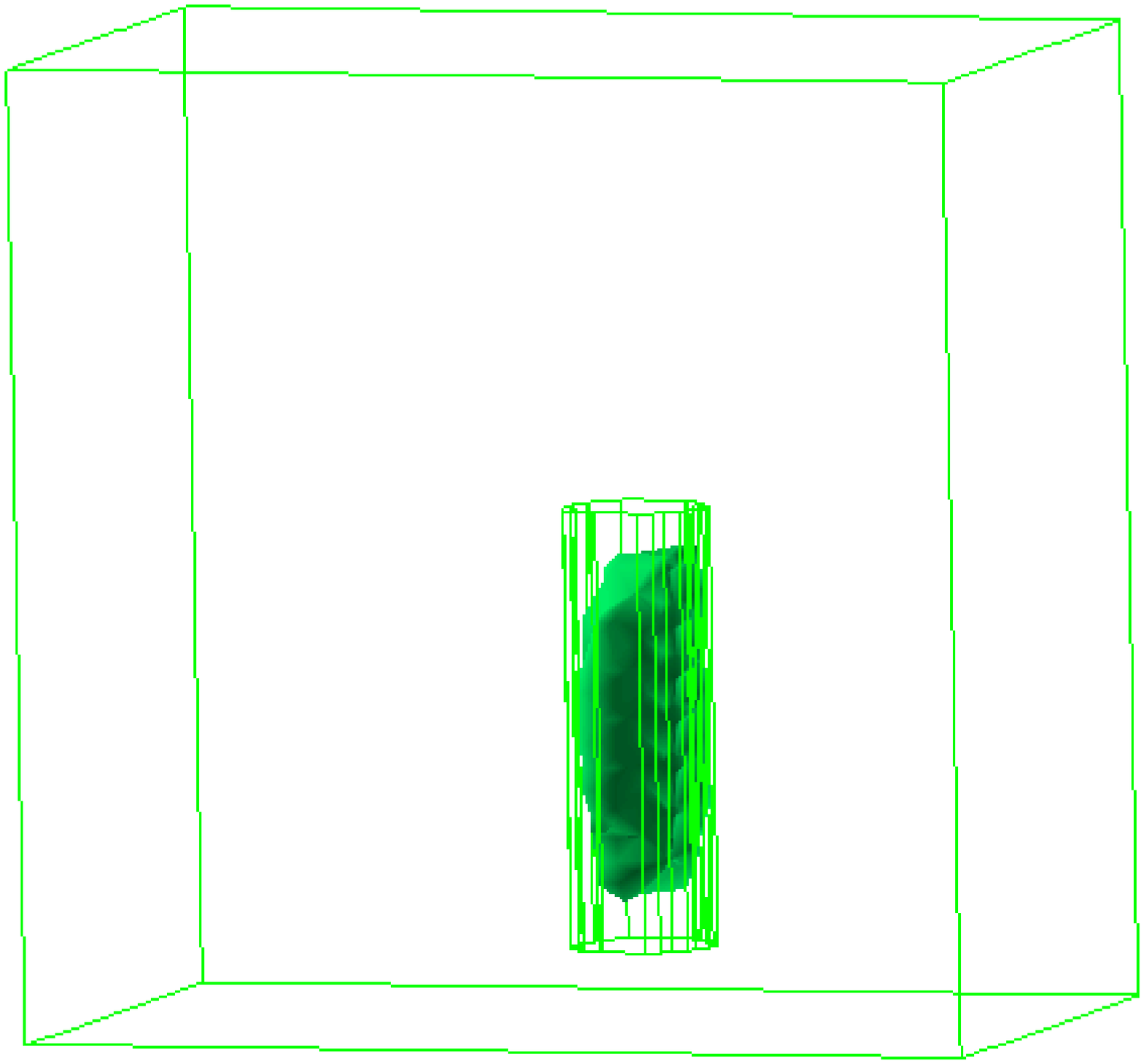}} &
{\includegraphics[scale=0.25,angle=-90, clip=]{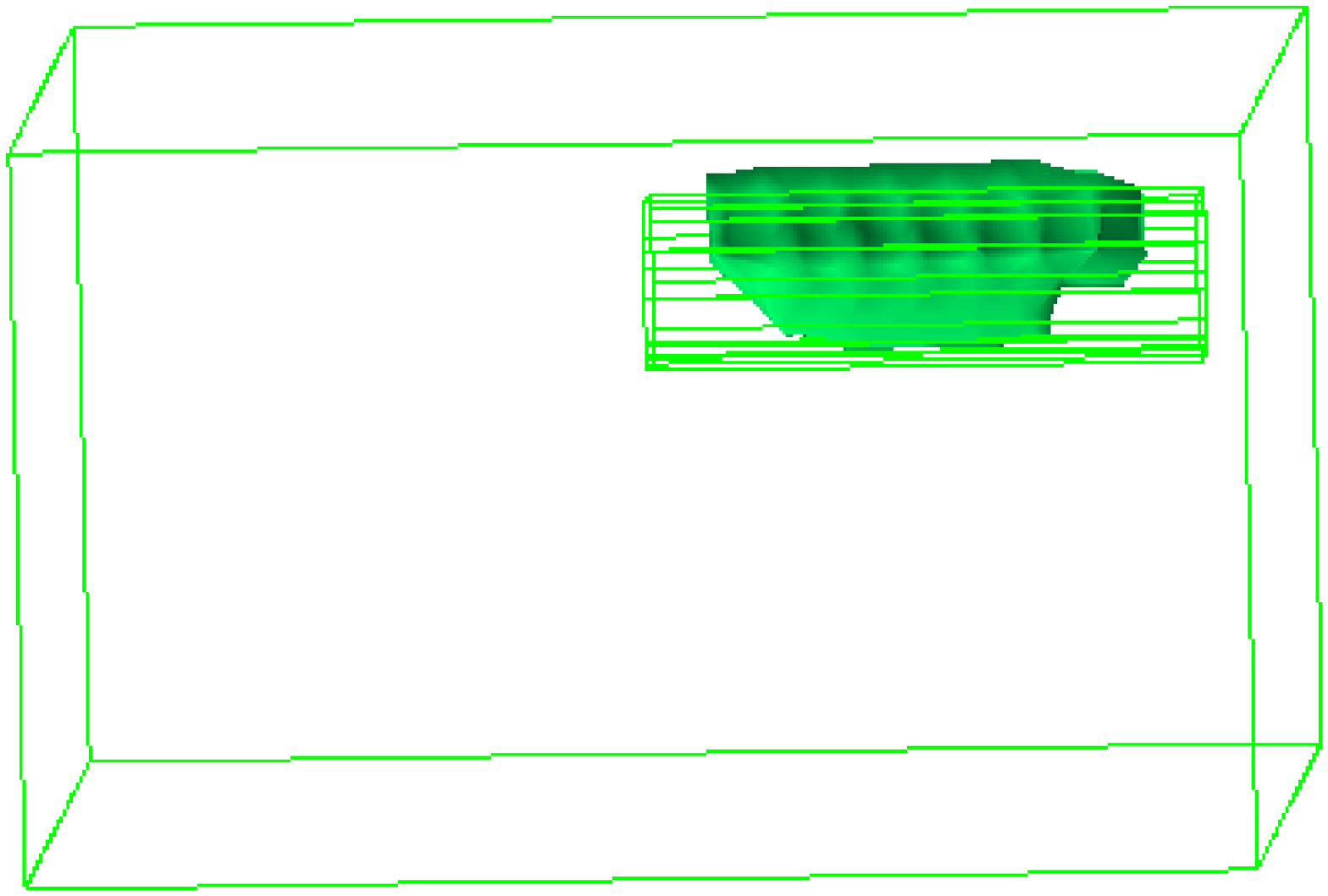}} \\
a) twice refined mesh, $xy$-view & b)  twice refined mesh, $yz$-view \\
\end{tabular}
\end{center}
\caption{{\protect\small \emph{Computed image of target number 2 of Table \ref{tab:table1}
.
Thin lines indicate correct
shape. To have a better visualization we have zoomed the domain $\Omega$ in
(\protect\ref{8.0}) in the domain $\Omega_\mathrm{FEM}$ in (\protect\ref{zoom}).
This target, which was a plastic bottle filled with water, was quite a large vertical size of 20 cm. On the other hand, our incident signal had a low power, which was much lower at the top and bottom of this target. This is why we were unable to image well the vertical size of this target. Still, one can observe that the image is stretched in the vertical direction.
}}}
\label{fig:3}
\end{figure}

\begin{figure}[tbp]
\begin{center}
\begin{tabular}{cc}
{\includegraphics[scale=0.19, angle=-90,clip=]{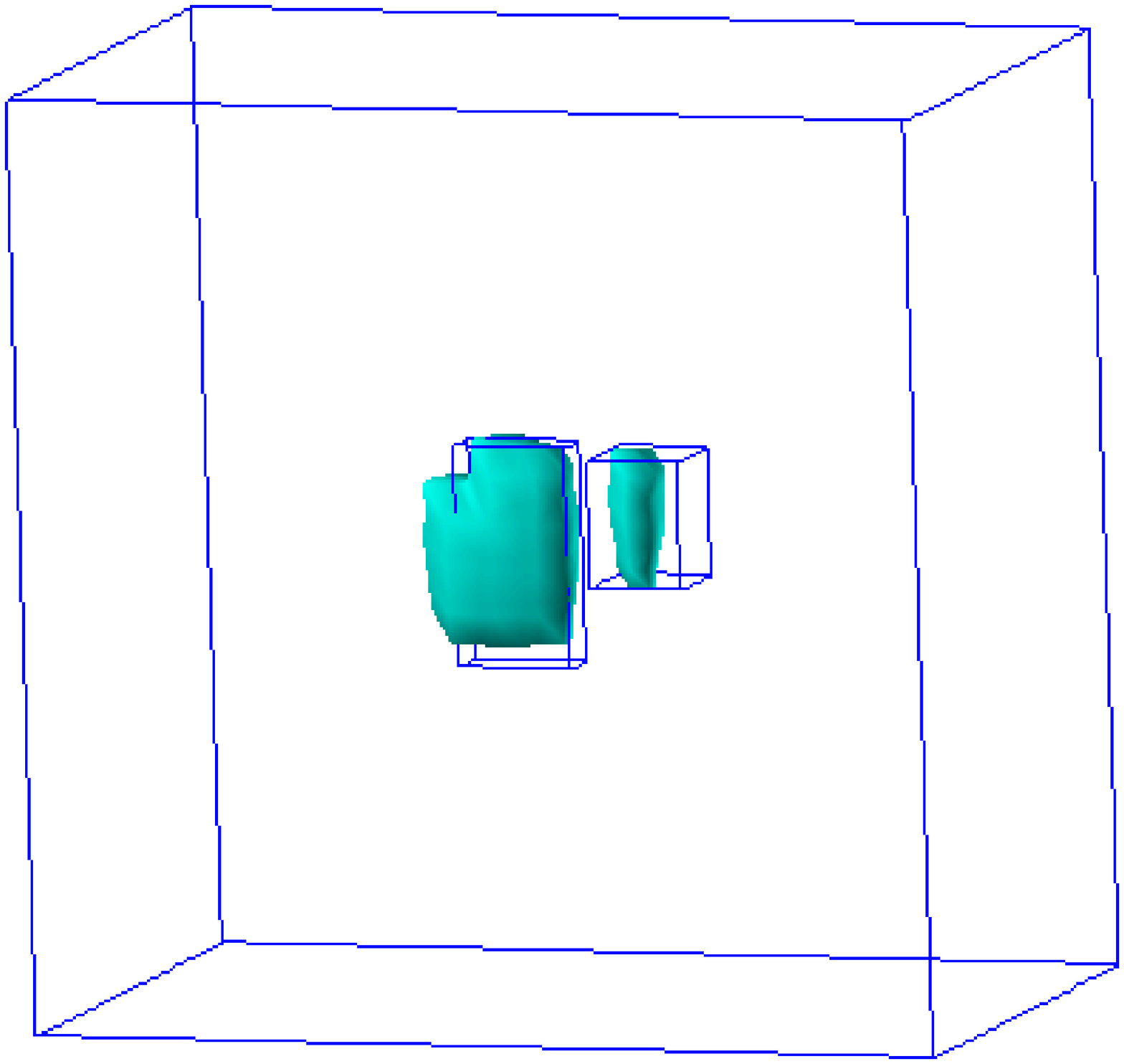}} &
{\includegraphics[scale=0.19, angle=-90,clip=]{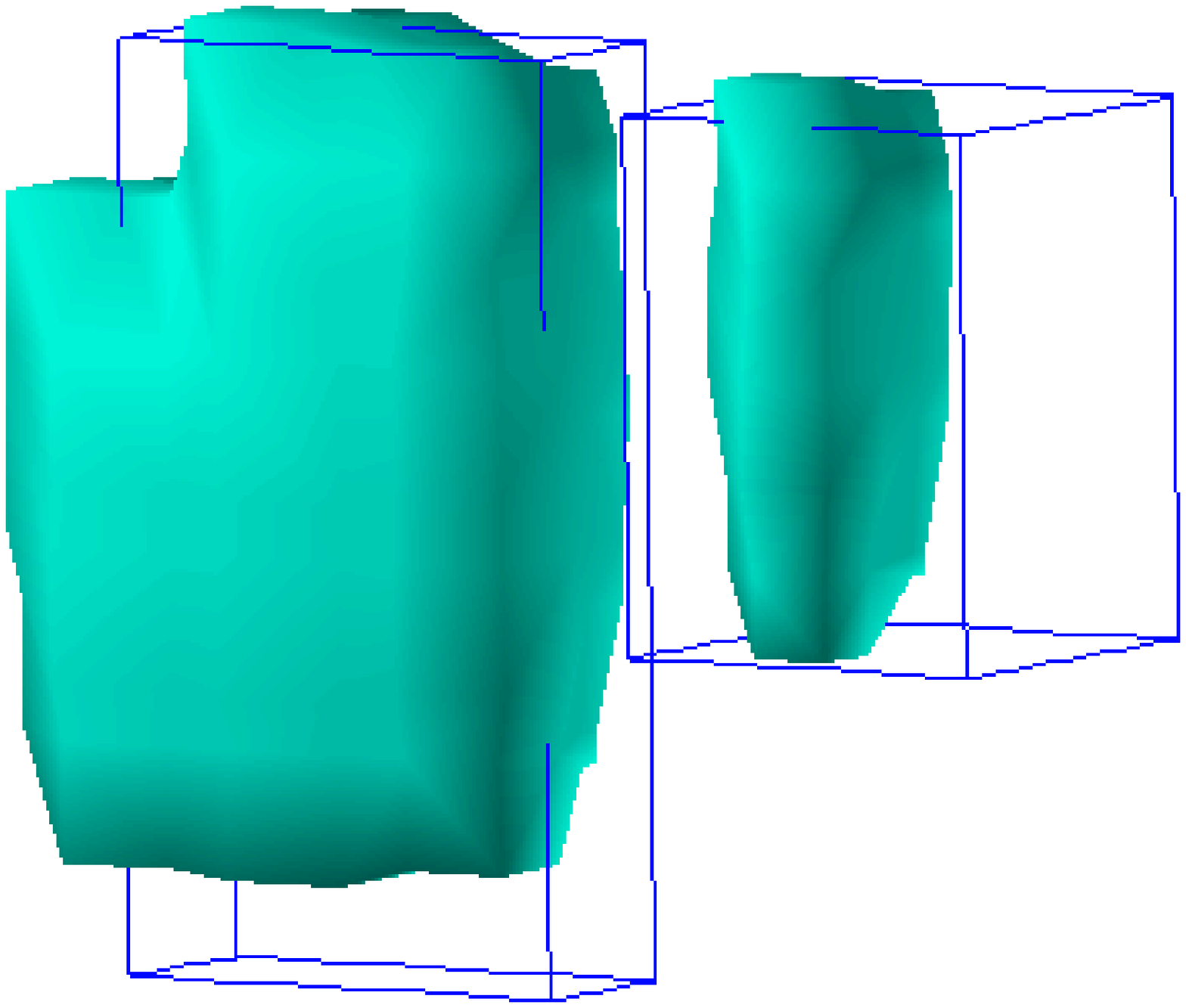}}\\
a) coarse mesh & b) zoomed  view \\
{\includegraphics[scale=0.19,angle=-90, clip=]{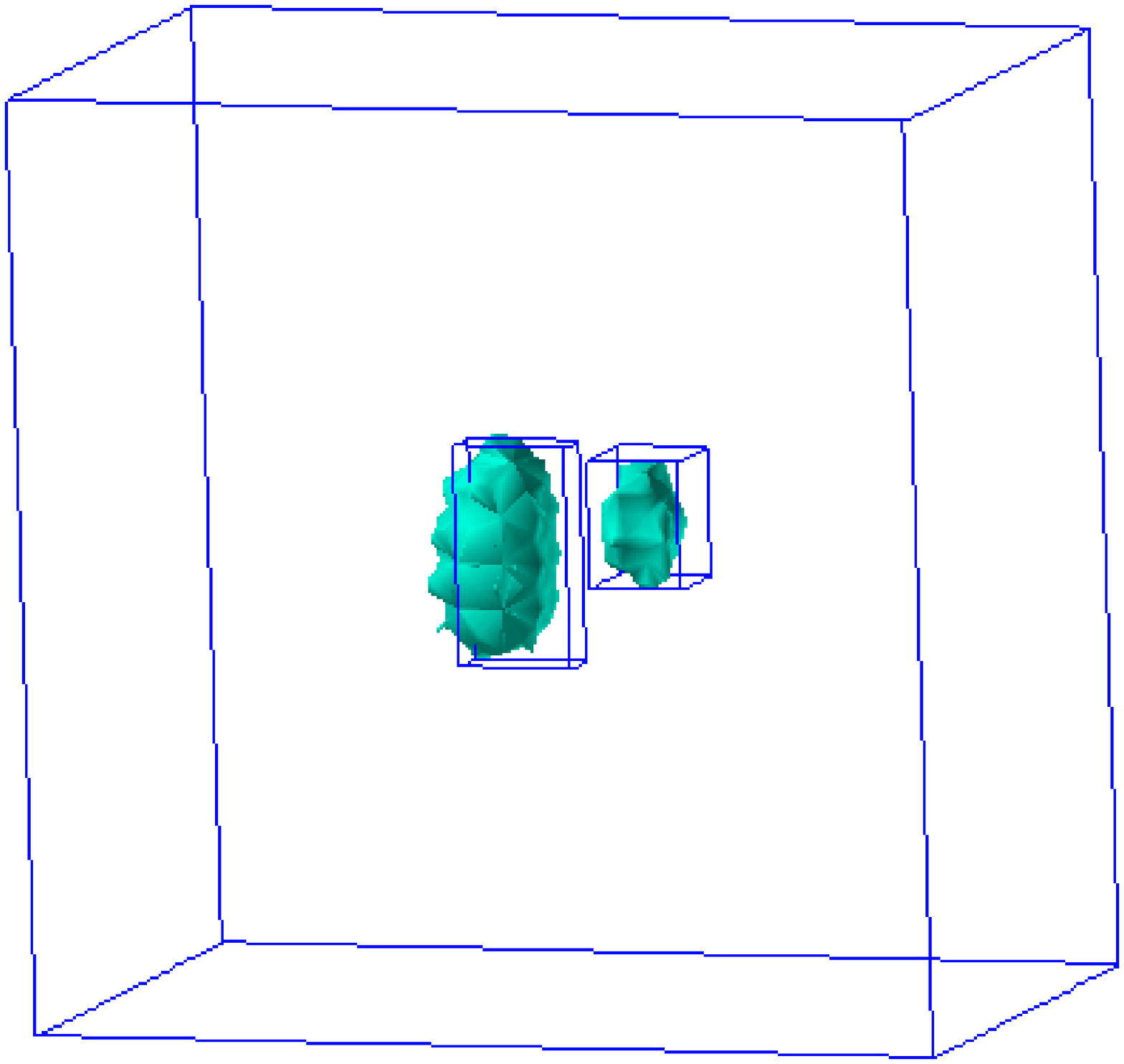}} &
{\includegraphics[scale=0.19,angle=-90, clip=]{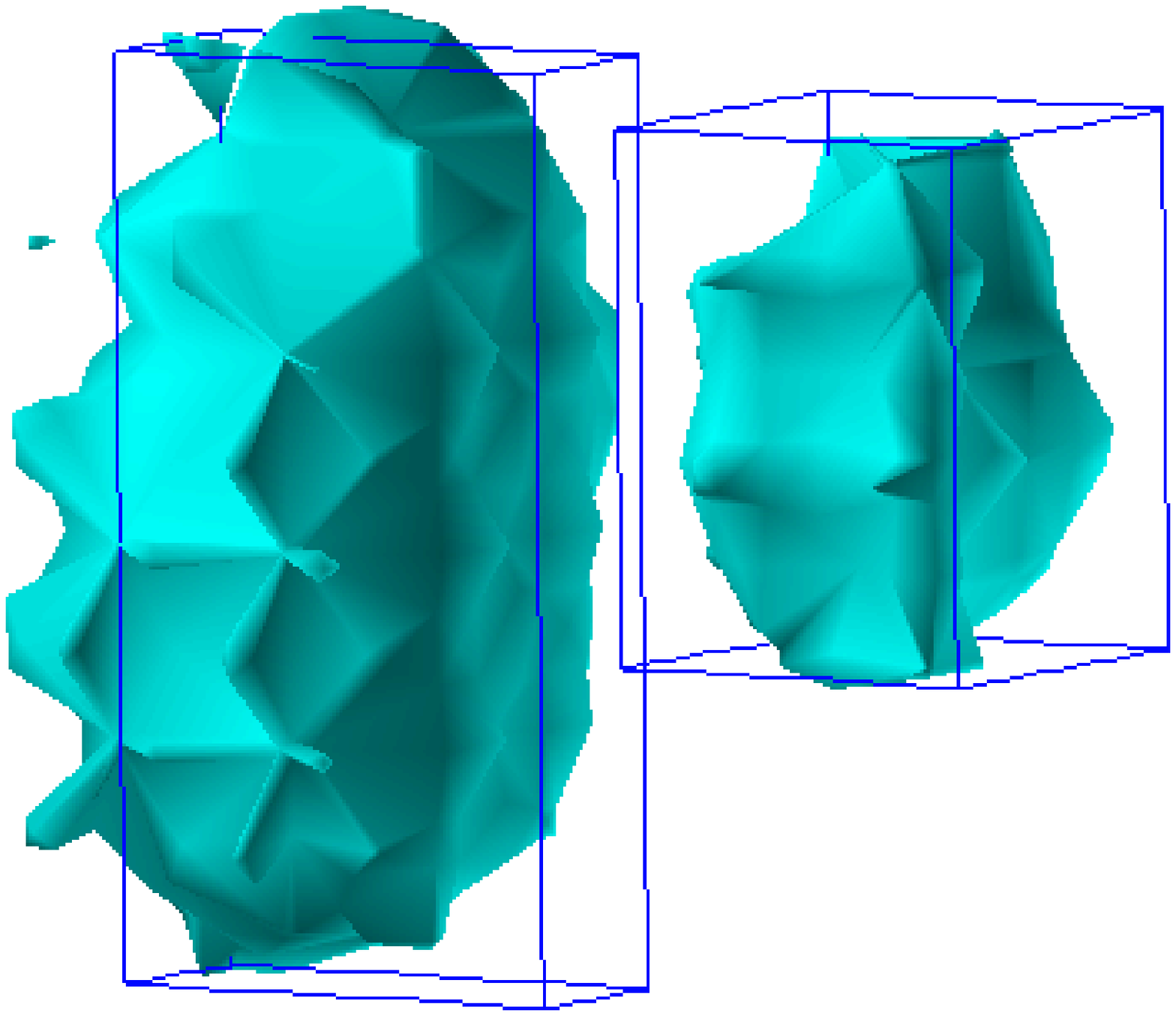}} \\
c)  two times refined mesh & d)  zoomed view \\
{\includegraphics[scale=0.19,angle=-90, clip=]{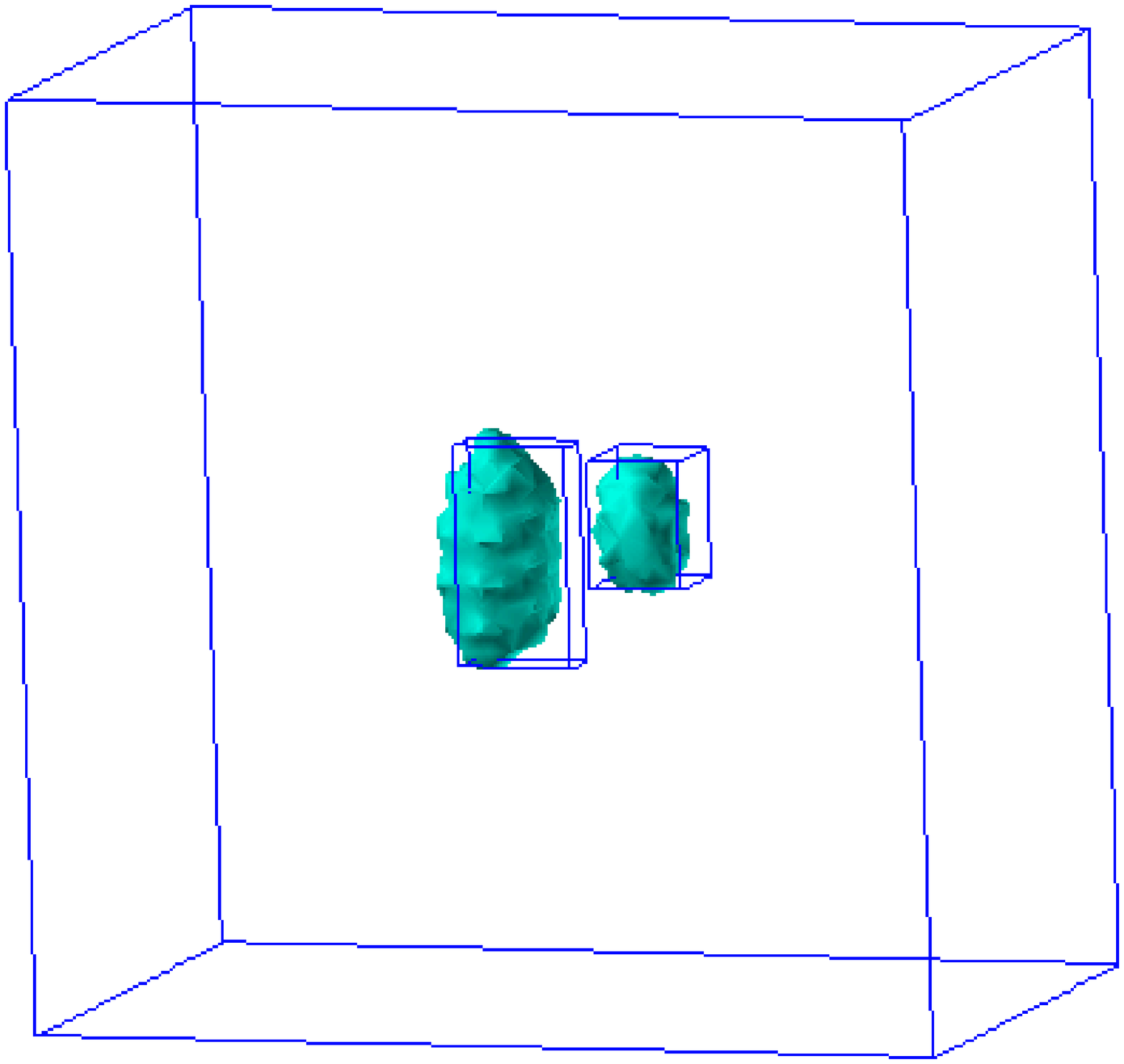}} &
{\includegraphics[scale=0.19,angle=-90, clip=]{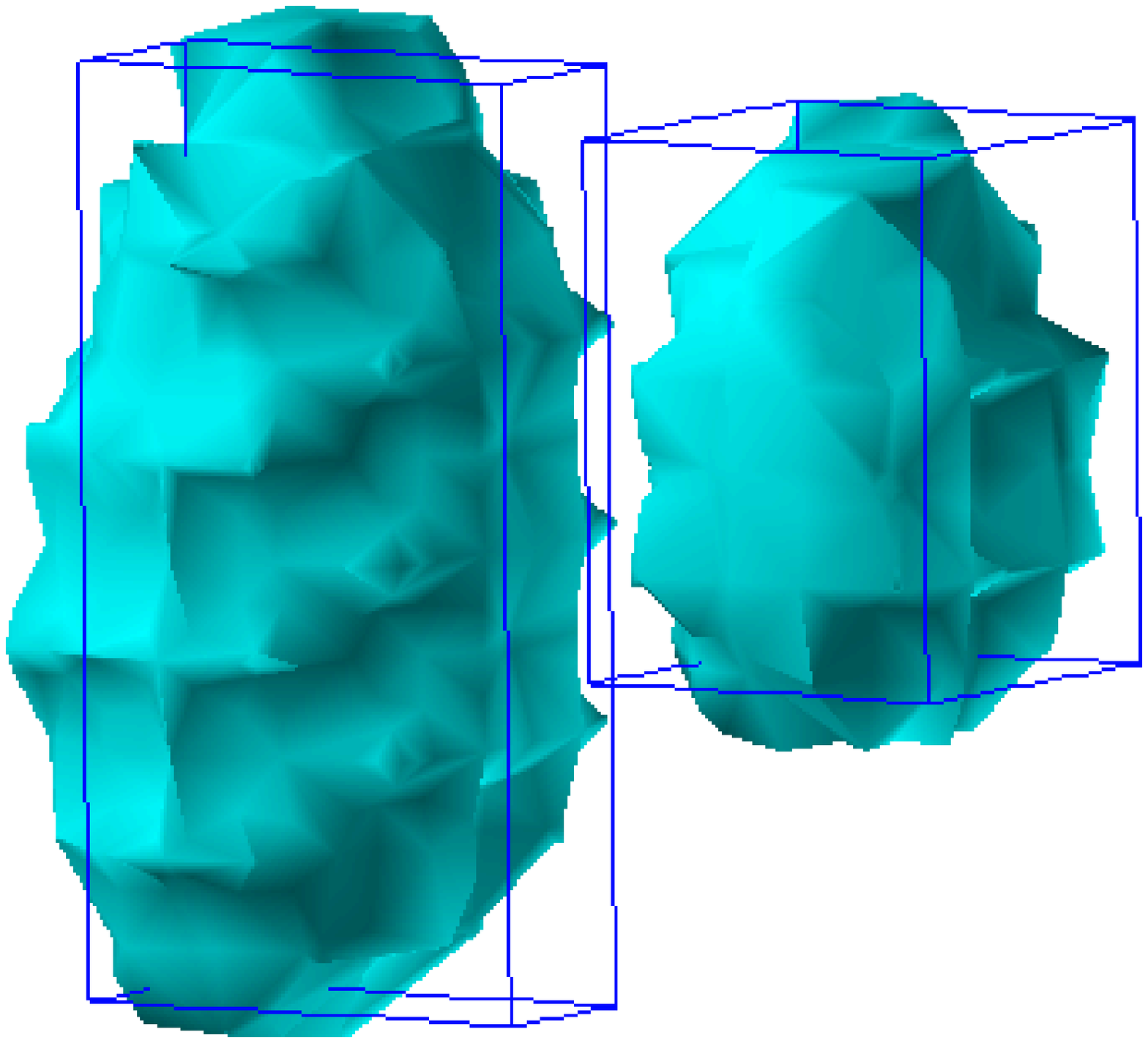}}  \\
e)  three times refined mesh  & f)  zoomed view \\
{\includegraphics[scale=0.19,angle=-90, clip=]{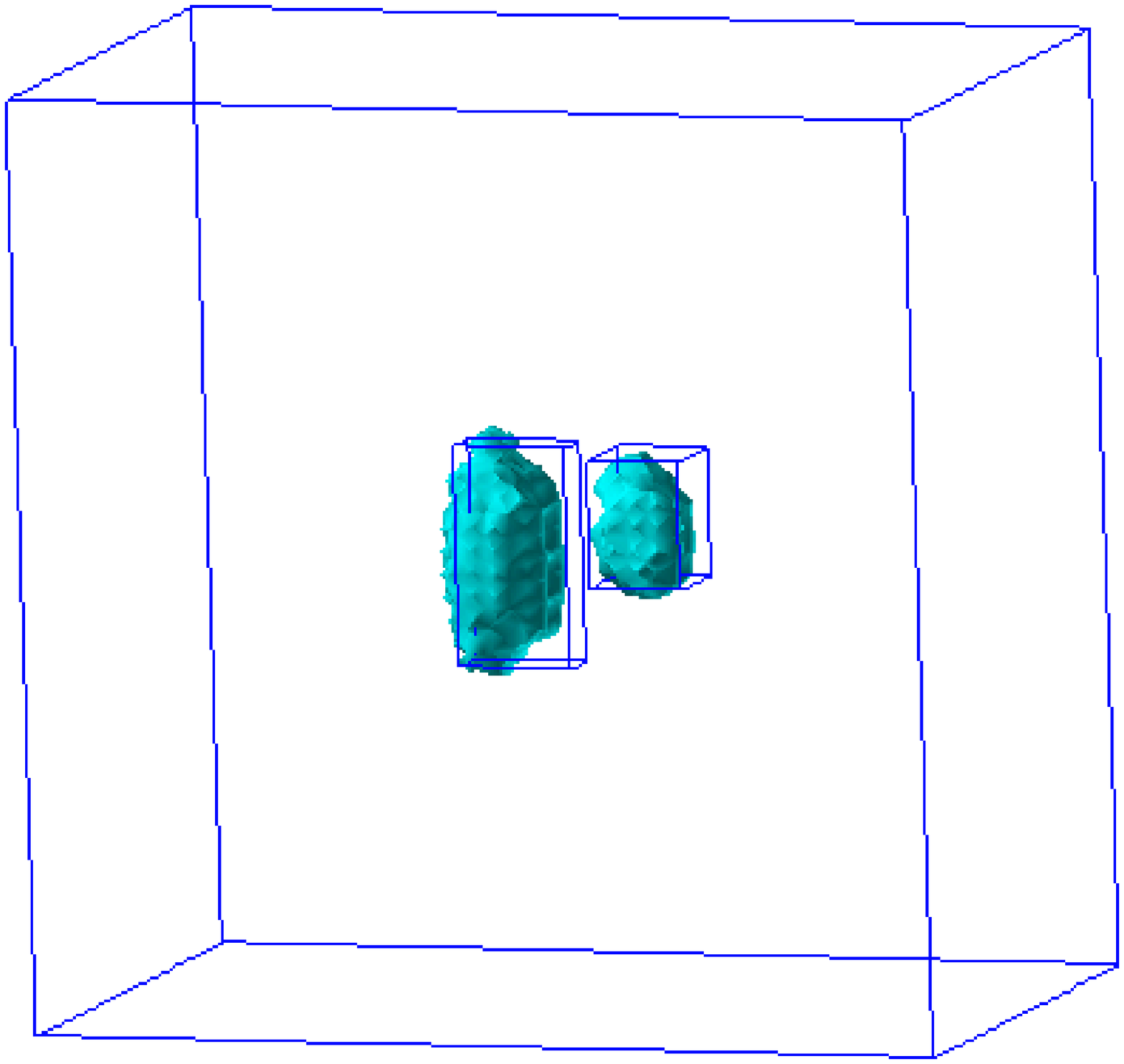}} &
{\includegraphics[scale=0.19,angle=-90, clip=]{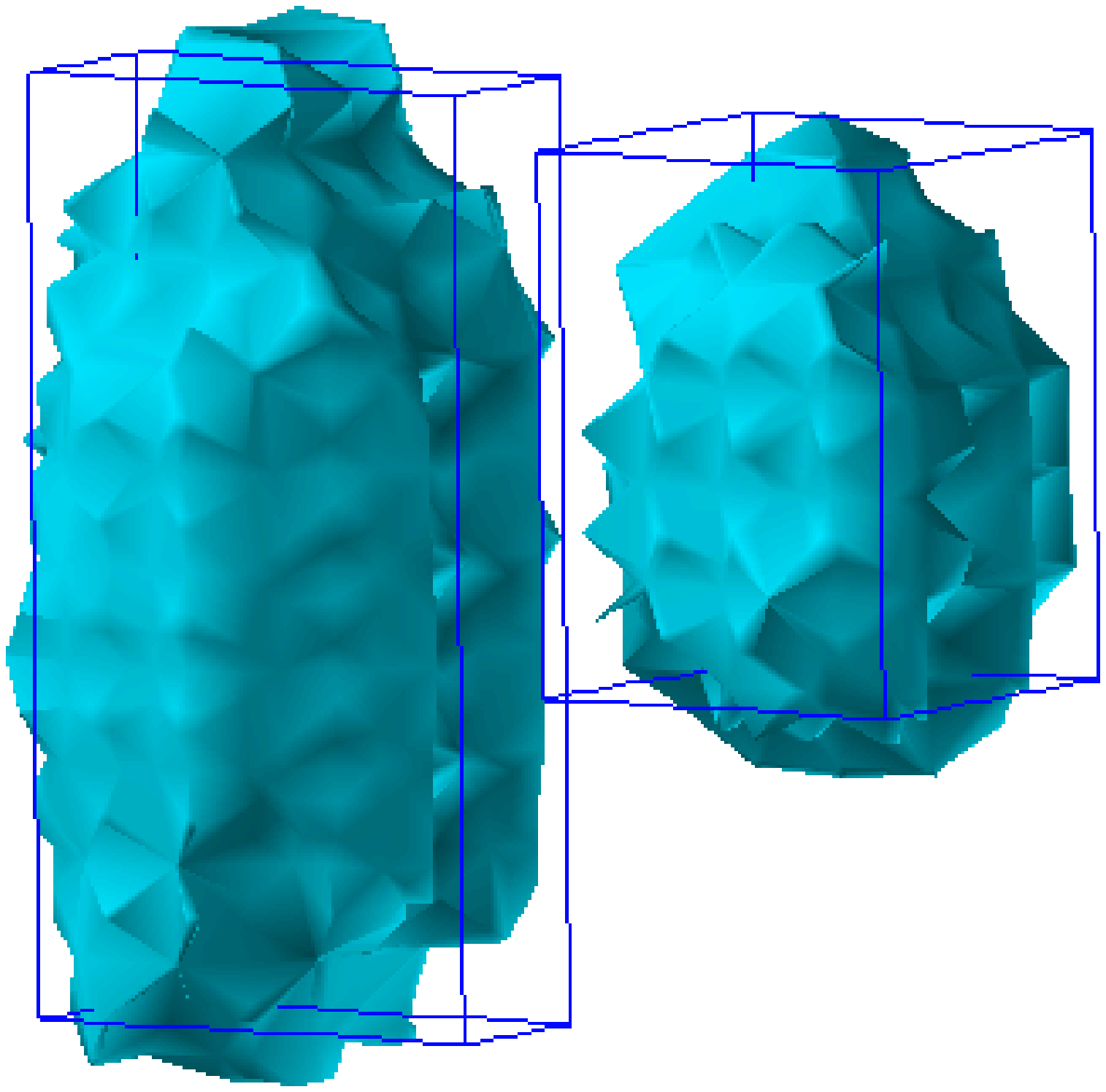}}
\\
g)  four times refined mesh & h) zoomed view
\end{tabular}
\end{center}
\caption{{\protect\small \emph{Computed images of targets number 4 of
      Table \ref{tab:table1} when superresolution is achieved on four
      times adaptively refined meshes.  Compare with
      Figure \ref{fig:1}-d).
}}}
\label{fig:5}
\end{figure}

Figures \ref{fig:2}--\ref{fig:7} display adaptively
refined meshes and 3D images of some targets of Table
\ref{tab:table1}.  To have a better visualization, we have zoomed some
figures from the domain $\Omega_\mathrm{FEM}$ defined in (\ref{8.0}) to the
domain defined in (\ref{zoom}).  We can conclude that the location of all
targets as well as their sizes in the $x$-, $y$-, and $z$-directions are well estimated
on the second stage of our two-stage numerical procedure.

\begin{table}[tbp]
\begin{center}
{\footnotesize \
\begin{tabular}{|c|l|l|l|}
\hline
Object & Blind/ & Description of target & Material \\
\# & Non-blind &  &  \\ \hline
1 & Non-blind & A metallic ball, 3 cm burial depth  & Metal \\ \hline
2 & Non-blind & A bottle filled with clear water, 3 cm depth & Water \\ \hline
3 & Blind & A ceramic mug, 5 cm burial depth  & Ceramic \\ \hline
4 & Non-blind & Two metallic blocks at 1 cm separation & Metal/Metal \\
\hline
\end{tabular}
}
\end{center}
\caption{{\protect\small \emph{Description of the data sets.}}}
\label{tab:table1}
\end{table}

\begin{table}[tbp]
\begin{center}
{\footnotesize \
\begin{tabular}{|c|l|r|r|l|l|l|}
\hline
Object & Material & Computed & Exact & Computed   & Measured
\\
\# &  & depth & depth & $n$  & $n$ \\
\hline
2 & Water & 3.6 & 4.0  & 4.7 & 4.88 \\ \hline
3 & Ceramic & 4.0 & 5.0  & 1.0 & 1.39 \\  \hline
\end{tabular}
}
\end{center}
\caption{{\protect\small \emph{Result of the first stage: the
      refractive indices $n = \protect\sqrt{\protect\varepsilon_\mathrm{r}}$ and the
      burial depths of non-metallic targets.}}}
\label{tab:table2}
\end{table}

\begin{table}[tbp]
\begin{center}
{\footnotesize \
\begin{tabular}{|c|l|r|r|l|l|}
\hline
Object & Material & Computed & Exact & Computed $\varepsilon_\mathrm{r}$ \\
\# &  & depth & depth  &  \\ \hline
1 & Metal & 2.9 & 3.0 & 31.0 \\ \hline
4 & Metal & 3.8 & 4.0 & 99.8 \\
& Metal & 4.0 & 4.0  & 56.5 \\ \hline
\end{tabular}
}
\end{center}
\caption{{\protect\small \emph{Result of the  first stage: the estimated effective dielectric constants and the burial depths of metallic targets. Object \#4 consists of two metallic
targets with 1 cm distance between their surfaces.}}}
\label{tab:table3}
\end{table}

\begin{table}[tbp]
\begin{center}
{\footnotesize \
\begin{tabular}{|l|l|l|l|l|l|c|}
\hline
Target number & 2 & 3   \\ \hline
blind (yes/no)  & no & no   \\ \hline
Measured $n$                      &  4.88  & 1.39   \\ \hline
$n^{\mathrm{comp}}$ coarse mesh      & 4.7  &  1  \\ \hline
$n^{\mathrm{comp}}$ 1 time ref. mesh  & 4.7  & 1  \\ \hline
$n^{\mathrm{comp}}$ 2 times ref.mesh   & 4.7  &  1 \\ \hline
$n^{\mathrm{comp}}$ 3 times ref.mesh   & 4.7  &  1 \\ \hline
\end{tabular}
}
\end{center}
\caption{\emph{ Stage 2. Computed $n^{\mathrm{comp}}$ and directly measured $n$ refractive
indices of dielectric targets.
 }}
\label{tab:table5}
\end{table}

\begin{table}[tbp]
\begin{center}
{\footnotesize \
\begin{tabular}{|l|l|l|l|l|l|l|l|}
\hline
Target number & 1 & 4   \\ \hline
blind (yes/no)  & no & no   \\ \hline
$\varepsilon_\mathrm{r}^{\mathrm{comp}}$ coarse mesh      & 24.5 & 75.6  \\ \hline
$\varepsilon_\mathrm{r}^{\mathrm{comp}}$ 1 time ref. mesh & 24.6 & 100  \\ \hline
$\varepsilon_\mathrm{r}^{\mathrm{comp}}$ 2 times ref.mesh & 24.7 & 100   \\ \hline
$\varepsilon_\mathrm{r}^{\mathrm{comp}}$ 3 times ref.mesh & 24.6 & 100  \\ \hline
\end{tabular}
}
\end{center}
\caption{\emph{Stage 2. Computed appearing dielectric constants }$\protect\varepsilon %
_\mathrm{r}^{\mathrm{comp}}$\emph{\ of metallic targets.  }}
\label{tab:table6}
\end{table}

\begin{figure}[tbp]
\begin{center}
    \begin{tabular}{c c c}
      \includegraphics[width = 0.33\textwidth, clip = true, trim = 3cm 1cm 3cm 1cm]{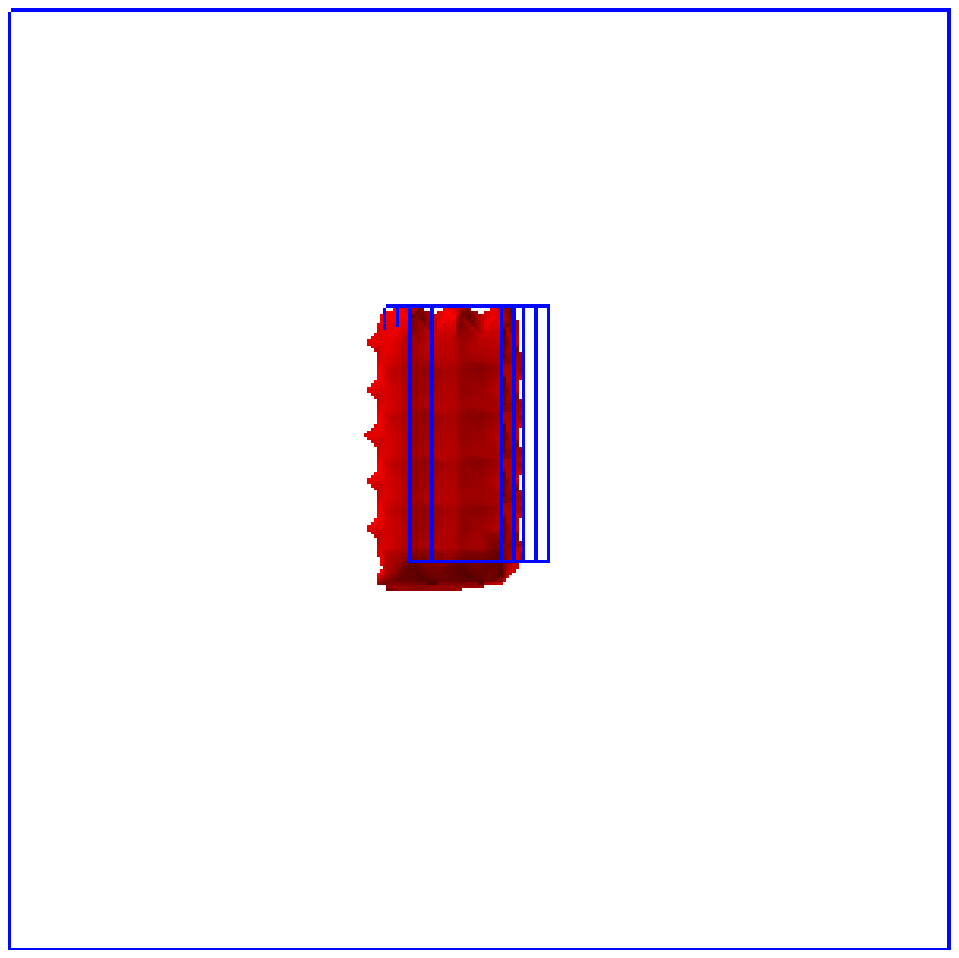} &
      \includegraphics[width = 0.33\textwidth, clip = true, trim = 3cm 1cm 3cm 1cm]{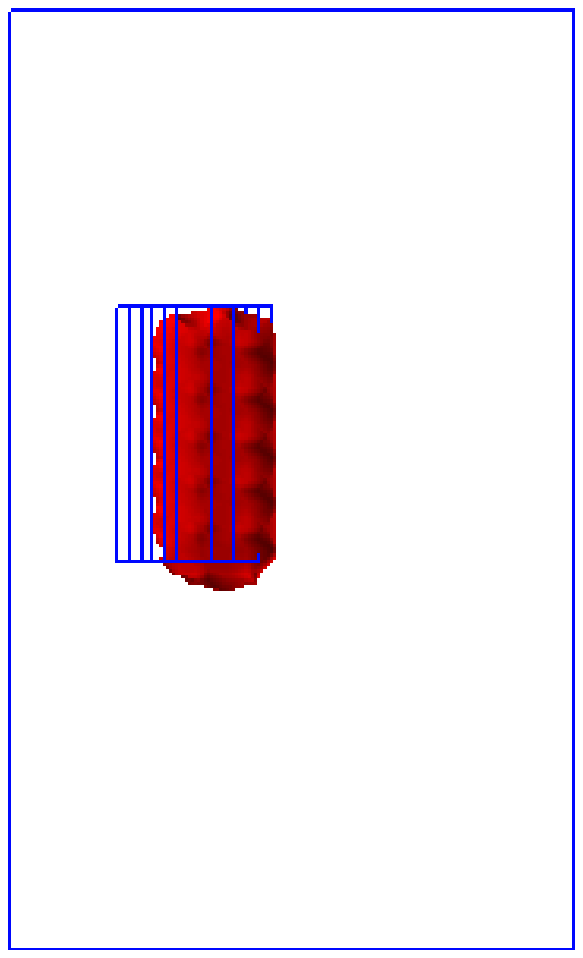} &
      \includegraphics[width = 0.33\textwidth, clip = true, trim = 2cm 0cm 2cm 0cm]{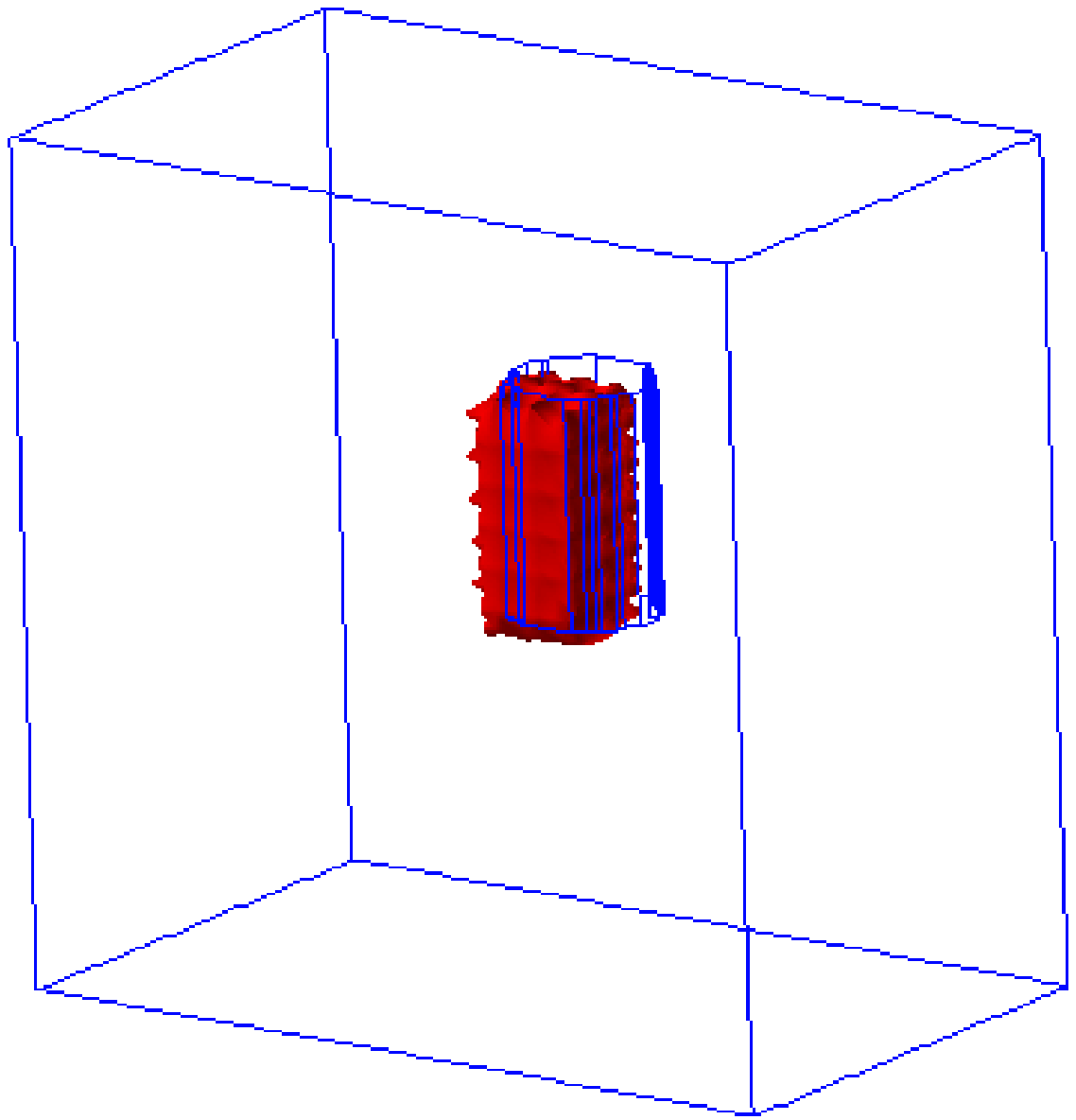} \\
      (a) Front & (b) Side & (c) Perspective \\
      \includegraphics[width = 0.33\textwidth, clip = true, trim = 4cm 3cm 4cm 1cm]{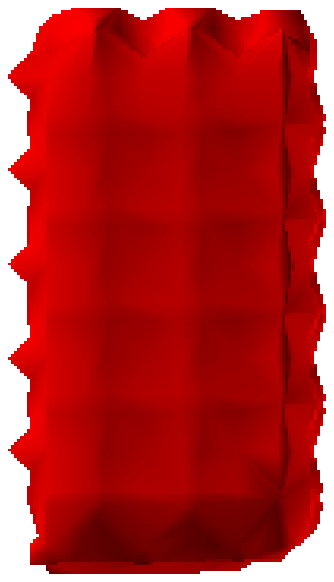} &
      \includegraphics[width = 0.33\textwidth, clip = true, trim = 3cm 3cm 5cm 1cm]{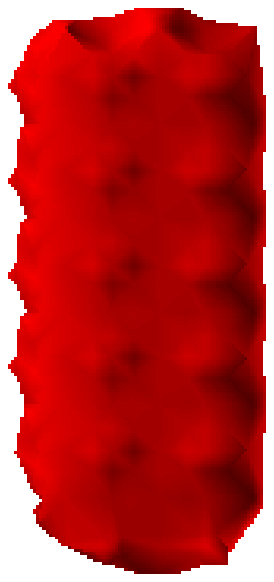} &
      \includegraphics[width = 0.33\textwidth, clip = true, trim = 4cm 3cm 3cm 0cm]{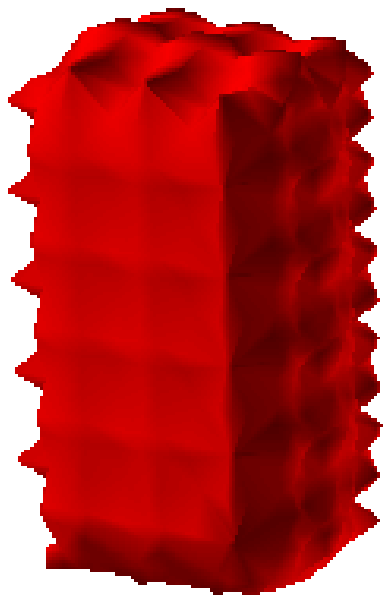} \\
      (d) Front, zoomed & (e) Side, zoomed & (f) Perspective, zoomed
    \end{tabular}
\end{center}
    \caption{Three views and zooms of the reconstruction of the target
      number 3 of Table \protect\ref{tab:table1} on the once refined
      mesh.  Recall that target number 3 is a  ceramic mug.}
  \end{figure}

\begin{figure}[tbp]
\begin{center}
    \begin{tabular}{c c c}
      \includegraphics[width = 0.33\textwidth, clip = true, trim = 3cm 1cm 3cm 1cm]{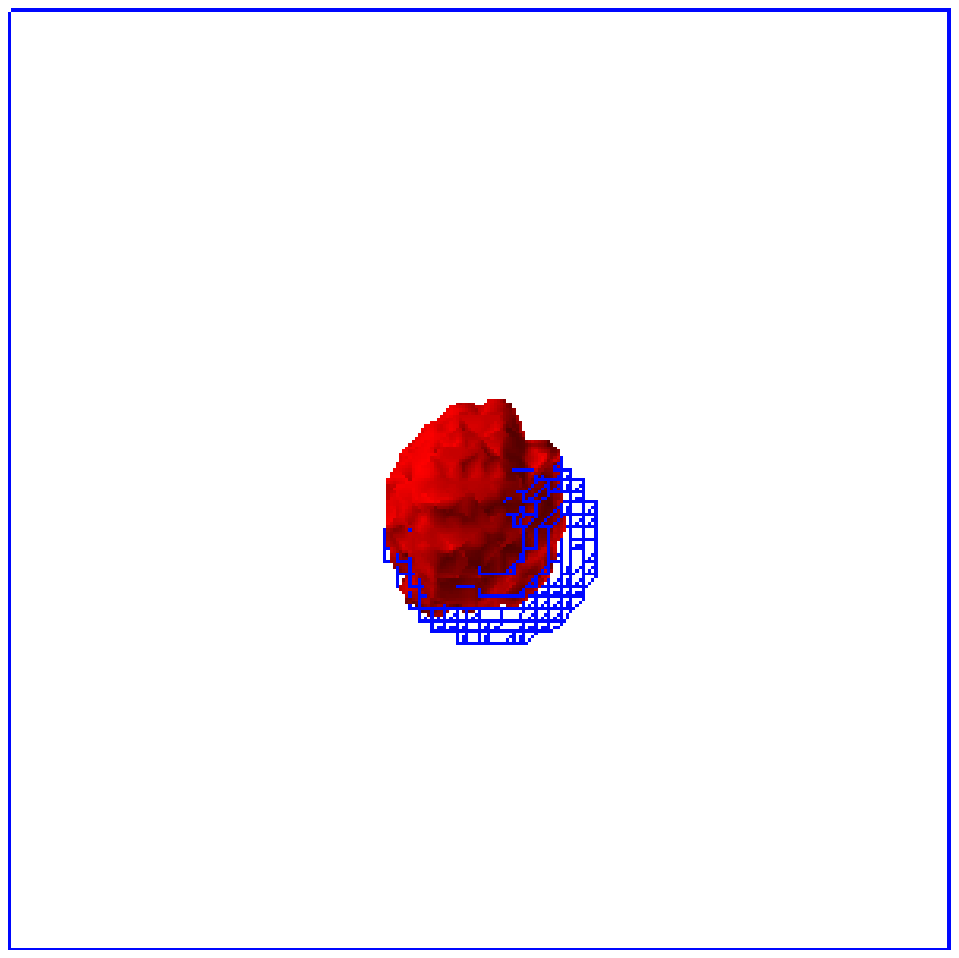} &
      \includegraphics[width = 0.33\textwidth, clip = true, trim = 3cm 1cm 3cm 1cm]{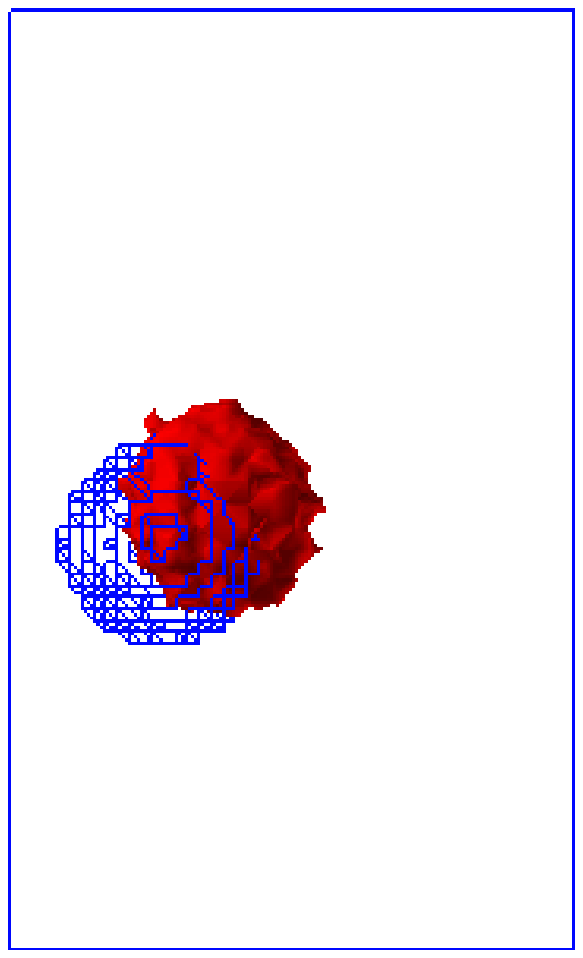} &
      \includegraphics[width = 0.33\textwidth, clip = true, trim = 2cm 0cm 2cm 0cm]{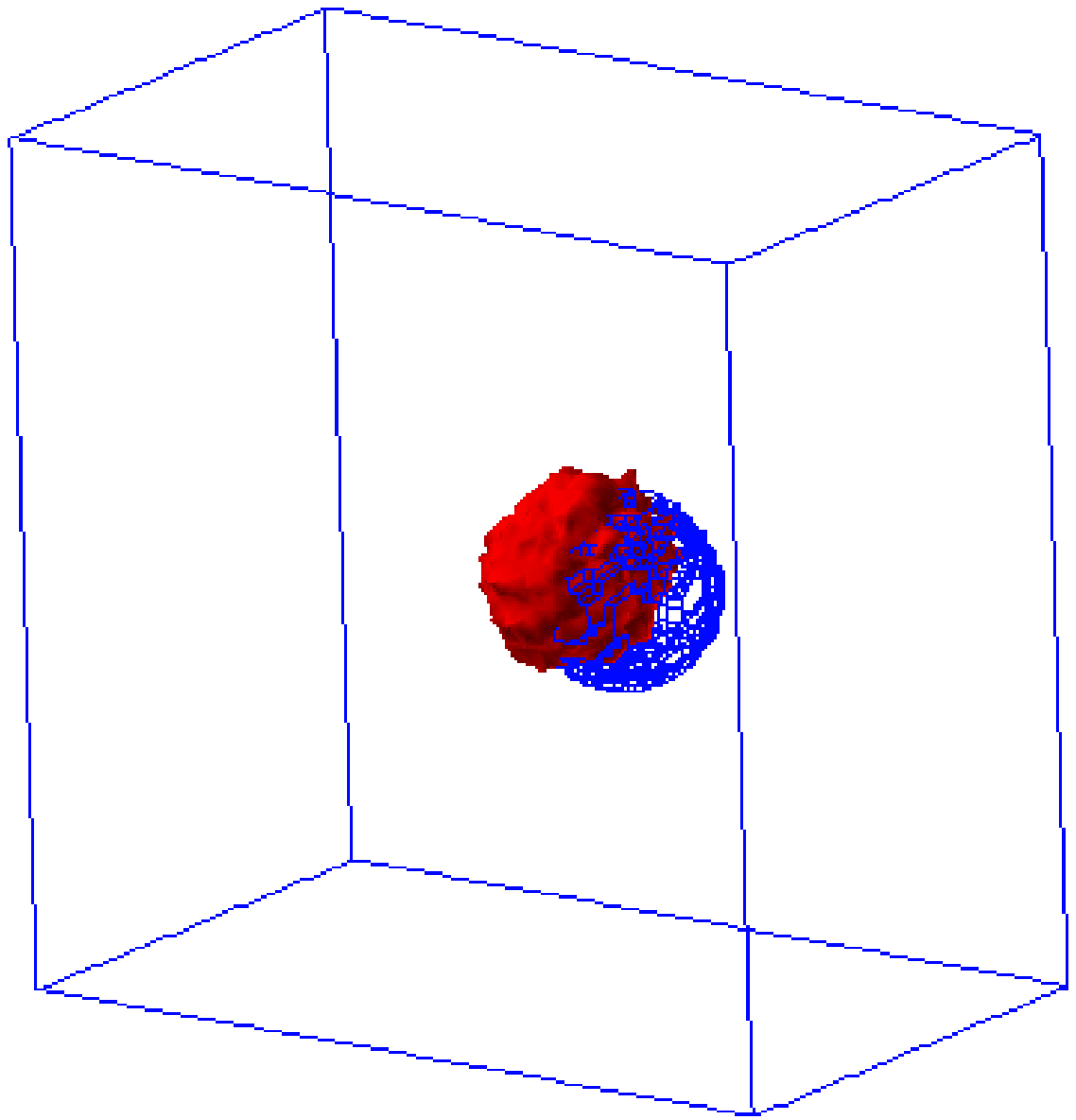} \\
      (a) Front & (b) Side & (c) Perspective \\
      \includegraphics[width = 0.33\textwidth, clip = true, trim = 4cm 2cm 4cm 2cm]{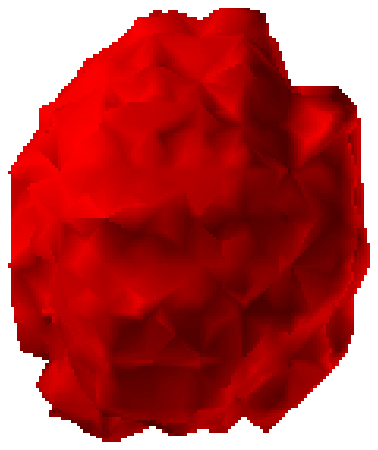} &
      \includegraphics[width = 0.33\textwidth, clip = true, trim = 3cm 2cm 5cm 2cm]{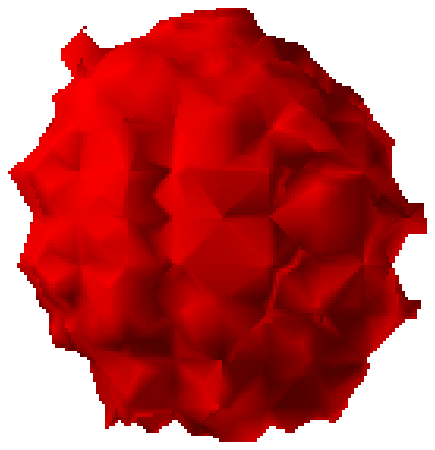} &
      \includegraphics[width = 0.33\textwidth, clip = true, trim = 4cm 2cm 3cm 1cm]{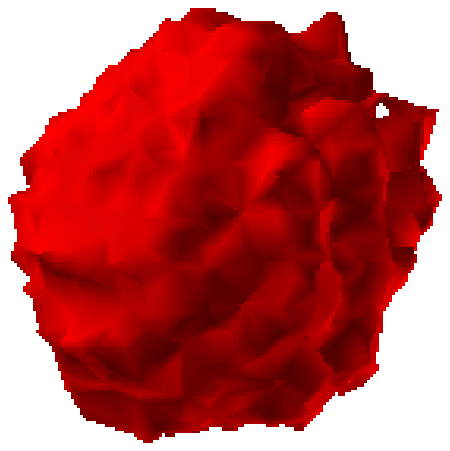} \\
      (d) Front, zoomed & (e) Side, zoomed & (f) Perspective, zoomed
    \end{tabular}
\end{center}
    \caption{Three views and zooms of the reconstruction of the target
      number 1 of Table \protect\ref{tab:table1} on the three times
      refined mesh. The initial guess in this test is taken from Test
      2 of \cite{NBKF2}, see Figure 5.1-b),d) of \cite{NBKF2}. Recall that target number 1 is a metallic ball.}
\label{fig:7}
  \end{figure}

\newpage

\section{Summary}

This is the fifth (5th) paper (after \cite{BTKF, BTKJ, NBKF, NBKF2}) in the recent
series of publications of this group about the performance of the
two-stage numerical procedure of \cite{BK} on experimental backscattering
time-dependent data generated by a single location of the source of
electromagnetic waves. While in \cite{BTKF, BTKJ, NBKF} we have considered the case
of targets placed in air, in \cite{NBKF2} and here we consider the more
challenging case of targets buried in the ground. This case is
more challenging because the signal scattered by the ground is
heavily mixed with the signal scattered by the target.

It was shown in \cite{NBKF2} that the globally convergent numerical method of
\cite{BK} accurately images refractive indices and locations of buried
targets. In this paper we complement the globally convergent method by
the locally convergent adaptivity technique. The adaptivity takes the
image of the globally convergent method as the starting point for
subsequent iterations. The theory of the adaptivity can be found in
\cite{AB, KBB, BK1, BK2, BKK, BJ1, B}, \cite{BMaxwell2, BJ, BK}, and \cite{BK_AA}. In particular, the important analytical
guarantee of the fact that adaptivity indeed refines images was first
established in [6] and then also published in \cite{BK} and \cite{BK_AA}.

As a result of the application of the adaptivity, our images are
significantly refined: the shapes of the targets are accurately imaged. A
particularly interesting case is the case of the superresolution
(Figure 2-d and Figures 5). We have accurately imaged both targets in
this case.

In conclusion, we believe that the two-stage numerical procedure of
\cite{BK} is now completely verified on experimental data.

\begin{center}
\textbf{Acknowledgments}
\end{center}

This research was supported by US Army Research Laboratory and US Army
Research Office grant
 W911NF-11-1-0399, the
Swedish Research Council, the Swedish Foundation for Strategic
Research (SSF) through the Gothenburg Mathematical Modelling Centre
(GMMC). The computations were performed on
resources at Chalmers Centre for Computational Science and Engineering
(C3SE) provided by the Swedish National Infrastructure for Computing
(SNIC).

\end{document}